\PassOptionsToPackage{numbers, compress}{natbib}
\documentclass{article}

\usepackage[utf8]{inputenc} 
\usepackage[T1]{fontenc}    
\usepackage{hyperref}       
\usepackage{color}
\usepackage{url}            
\usepackage{booktabs,multirow}       
\usepackage{amsmath,amssymb,amsfonts}       
\usepackage{mathrsfs}
\usepackage{amsthm}
\usepackage{nicefrac}       
\usepackage{microtype}      
\usepackage{listings}
\usepackage{tikz}
\usepackage[ruled,vlined]{algorithm2e}
\usepackage{graphicx,float}
\usepackage{geometry}
\usepackage{bm}
\usepackage{indentfirst}

\DeclareFontFamily{U}{rsfs}{\skewchar\font127 }
\DeclareFontShape{U}{rsfs}{m}{n}{%
   <-6> rsfs5
   <6-8> rsfs7
   <8-> rsfs10
}{}
\renewcommand{\d}{\mathrm{d}}

\newcommand{\pt}{\partial}

\newcommand{\e}{\mathrm{e}}

\newcommand{\<}{\left<}
\renewcommand{\>}{\right>}

\newcommand{\ran}{\mathcal{R}}
\newcommand{\mcalP}{\mathcal{P}}
\newcommand{\BbbP}{\mathbb{P}}
\newcommand{\BbbR}{\mathbb{R}}

\newcommand{\intbar}{\rlap{\,---}\int}
\newcommand{\smallintbar}{\rlap{\,\(-\)}\int}

\newtheorem{theorem}{Theorem}
\newtheorem{lemma}[theorem]{Lemma}
\newtheorem{proposition}[theorem]{Proposition}
\newtheorem{corollary}[theorem]{Corollary}
\theoremstyle{definition}
\newtheorem{remark}{Remark}
\newtheorem{definition}{Definition}
\newtheorem{example}{Example}

\numberwithin{equation}{section}

\begin{document}

\title{A convergent numerical algorithm for the stochastic growth-fragmentation problem}
\author{Dawei Wu
\footnote{School of Mathematical Sciences, Peking University, Beijing, 100871, P.R. China. Email: \url{2201110052@pku.edu.cn}} \and
Zhennan Zhou
\footnote{Beijing International Center for Mathematical Research, Peking University, Beijing, 100871, P.R. China. Email: \url{zhennan@bicmr.pku.edu.cn}}}
\date{Dec 18, 2022}
\maketitle

\begin{abstract}
    The stochastic growth-fragmentation model desribes the temporal evolution of a structured cell population through a discrete-time and continuous-state Markov chain.
    The simulations of this stochastic process and its invariant measure are of interest.
    In this paper, we propose a numerical scheme for both the simulation of the process and the computation of the invariant measure, and show that under appropriate assumptions, the numerical chain converges to the continuous growth-fragmentation chain with an explicit error bound.
    With a triangle inequality argument, we are also able to quantitatively estimate the distance between the invariant measures of these two Markov chains.
\end{abstract}



\section{Introduction}
\subsection{Overview of the growth-fragmentation model}
Growth-fragmentation models (abbreviated as ``GF'' henceforth) describe the temporal evolution of a structured population characterised by state variables such as age, size, etc..
As its name suggests, it integrates a ``growth process'' and a ``fragmentation process''. The state variables evolves according to some deterministic equation, but are subject to sudden changes (usually into parts) occurring at stochastic moments. The GF model appears in physics to describe the degradation phenomenon in polymers, droplets and in telecommunications systems to describe some internet protocols \cite{perthame2005}.

We will focus the GF model on cell division. As is supported by experimental evidence in \cite{robert2014}, cells in a culture medium increase in size deterministically and split into two offsprings randomly. The rates at which they grow or break are both determined by their current sizes (state variable).
A video from the supplementary material of \cite{wang2010} illustrates this process.
A common assumptions is that at each fragmentation, the cell splits into two offsprings with equal sizes \cite{doumic2009}. Scientists have also analyzed the case where the sizes of the offsprings are random and determined by a kernel \cite{doumic2010}.

\subsection{Stochastic growth-fragmentation model}
There are two main approaches to modelling the GF process. One is through a evolutional partial differential equation regarding the size distribution of the population \cite{perthame2005,perthame2007,doumic2010}; the other is through the genealogical tree which is Markovian \cite{doumic2015,hoffmann2016,bertoin2017,bernard2019}. We will be using the Markovian model proposed in \cite{doumic2015} in this paper.
Their motivation was the availability of observation schemes at the level of cell individuals, for example, the experiments in \cite{wang2010} where scientists prepared a thin vertical tube with the top open and the bottom sealed, fixed an E.coli cell at the sealed end and watched it grow. When the cell split, exactly one of its descendents was kept at the bottom of the tube, and thus we were able to track a single lineage of cells through time by recording a video of the whole process and identifying those bottom cells.
In fact, the new observation scheme did improve on previous rates of convergence obtained by indirect measurements, as in \cite{doumic2012}.
We will state a slightly different version of their model in this subsection.

Use an infinite binary tree $\mathbb{T}=\bigcup_{n=0}^\infty \{0,1\}^n$ to represent the family of cells. Each node in $\mathbb{T}$ is a finite array of 0's and 1's, and its two children append a 0 and 1 to this array respectively. For example, the node $01$ has two child nodes, $010$ and $011$.

Assign a set of variables $(\xi_u, \zeta_u, b_u)$ to each cell $u\in\mathbb{T}$.  $\xi_u$ stands for its size at birth, $\zeta_u$ stands for its total life span until the fragmentation, and $b_u$ is the moment when it's born. Its size as a function of time is denoted by $x_u(t), t\ge b_u$. A cell grows according to a deterministic rate $g(x)$, and at a certain time it breaks into two equal-sized offsprings. Therefore, if $u^+$ is either of the daughter cells of $u$, we have the basic relation
\begin{equation}
    2\xi_{u^+}=x_u(b_u+\zeta_u).
\end{equation}
The stochasticity of this model comes from the randomness of $\zeta_u$ when $\xi_u$ is given.

The size $x_u$ of a cell $u$ is a function of time  that satisfies an ordinary differential equation.
\begin{equation} \label{grow-ode}
    \begin{cases}
        \dot x_u=g(x_u), \\
        x_u(b_u)=\xi_u.
    \end{cases}
\end{equation}
The function $g(x)$ is termed the \textbf{growth rate}.
Its solution is given by the inverse relation
$$t=\int_{\xi_u}^x \frac{\d y}{g(y)}.$$

During each infinitesimal time range $[t,t+\d t]$ it breaks into two equal parts with probability $B(x_u)\d t$. That is,
$$\d \mathbb{P}(\zeta_u=t|\zeta_u\ge t, \xi_u)
=B(x_u(t)) \d t.$$
The function $B(x)$ is termed the \textbf{fragmentation rate}. Using the solution to \eqref{grow-ode} and some elementary calculus, we can compute the complementary cumulative distribution of the terminal size $x_u(\zeta_u)$.
\begin{equation}
    \mathbb{P}(x_u(\zeta_u)>x|\xi_u)=\exp\left( -\int_{\xi_u}^x \frac{B(s)}{g(s)}\d s \right).
\end{equation}

Therefore, using the relation $2\xi_{u^+}=x_u(b_u+\zeta_u)$, both being the size of the cell $u$ at fragmentation, we know that the distribution of the initial size of $u^+$ satisfies
\begin{equation} \label{trans-tail-0}
    \BbbP(\xi_{u^+}>y| \xi_u=x)=\exp\left( -\int_x^{2y} \frac{B(s)}{g(s)}\d s \right), \quad \forall y>\frac{x}{2}.
\end{equation}
This distribution is supported on $[\frac{\xi_u}{2},\infty)$.
Thus the initial sizes $\{\xi_u\}_{u\in\mathbb{T}}$ constitute a Markovian binary tree structure, where a child node depends only on the parent node.
However, we only keep track of one descendent at each splitting of cell so that only one lineage is observed, and the model becomes a Markov chain, which is simpler to analyze.

Denote by $\xi^n$ the $n$-th member along any lineage, so the dependence between $\xi^{n}=x$ and $\xi^{n+1}=y$ is specified by the relation \eqref{trans-tail-0}.
Therefore, the density function of the conditional distribution $\xi^{n+1}|\xi^n$.
\begin{equation} \label{trans-density}
    p(x,y)=I_{y\ge \frac{x}{2}} \frac{B(2y)}{g(2y)}\exp\left( -\int_x^{2y} \frac{B(s)}{g(s)}\d s \right).
\end{equation}
Then, we can get the Markovian kernel
\begin{equation} \label{trans}
    \mcalP(x,A)=\int_A p(x,y)\d y.
\end{equation}
Later on, this Markov chain model will be referred to as the \textbf{growth fragmentation chain}. Our version of the model differs from that in \cite{doumic2015} in the following two aspects: (1) our growth rate is an arbitrary function $g(x)$ while theirs is a linear function $\tau x$; (2) our growth rate is deterministic for all cells while their rate $\tau$ changes randomly when transitioning to an offspring.
\cite{doumic2015} has proved the $V$-uniform ergodicity \cite[Chapter 16]{mt1993} of their GF chain under certain conditions, which we will also review in Section 3.

The main topic of interest in the study of GF models lies in the inverse inference of $B(x)$ upon function values of $g(x)$ (which is usually known)) and the observation of the cell population.
\cite{perthame2007,doumic2009,doumic2012,doumic2015} have proposed inference schemes and proved the convergence rates.
In particular, \cite{doumic2015}, with its Markovian model, managed to improve the convergence rate by a slight order of magnitude in $N$ (the number of samples observed), showing the worth of the Markov chain model in solving the GF problem.
The authors thereof numerically tested the convergence through computer-generated samples, but did not carry out a rigorous analysis of the correctness of the numerical scheme.
Also, some approaches to the inverse problem, e.g. Bayesian statistics and the maximum likelihood estimate \cite{stuart2010,inigo2020}, still requires accurate computation of the forward problem of generating samples.

This situation poses a need for an efficient generator of the GF chain and its invariant measure. The numerical scheme we are about to propose in this paper is a finite volume method in discrete time, discretizing general probability distributions into piecewise constant measures on a truncated domain, which are essentially finite-dimensional. 
The numerical method we are to propose contains the effects of both the truncation of the domain onto $[0,a]$ and the discretization over a grid with mesh size $h$ in the proof of the convergence rate. The method is able to approximate the invariant measure of the GF chain with a high precision. 
More precisely, if we denote by $\pi$ the invariant measure of the GF chain and $\pi_{a,h}$ the invariant measure generated through our numerical scheme, we will prove that
$$\|\pi-\pi_{a,h}\|=O(h+ah+\e^{-a^{\alpha}})$$
under appropriate assumptions on $B(x)$ and $g(x)$.
In addition, the method with its truncated domain preserves the probability measure naturally and avoids using artificial boundary conditions usually required in differential equation based approaches. 

\subsection{Organization of this paper}
In section 2, we will discuss the discretization scheme in detail. It involves setting up a discrete state space and defining discrete transition kernels therein.
In section 3, we will introduce some basic results on the ergodicity of our model, and propose the main theorem about the convergence of the scheme.
In section 4, we will present the most crucial part of the proof.
In section 5, we will carry out numerical experiments to verify the convergence properties.

\subsection{Notations and symbols}
Here and throughout the whole paper, $\BbbR_+=(0,\infty)$ denote the positive real line as our state space. The probability measures on it is denoted by $\mathscr{P}(\BbbR_+)$. 
The distribution law of a random variable $\eta$ as a probability measure over $\BbbR_+$ will be denoted by $\mathcal{L}(\eta)$, and we also say that $\eta\sim \mathcal{L}(\eta)$.

The range of a linear operator $T$ is denoted by $\ran(T)$.

A Markovian kernel $\mcalP(x,A)$ over a measurable space $(\mathsf{X},\mathfrak{F})$ is a bivariate function over $\mathsf{X}\times \mathfrak{F}$ satisfying
\begin{enumerate}
    \item $\mcalP(x,\cdot)$ as a function on $\mathfrak{F}$ is a probability measure.
    \item $\mcalP(\cdot, A)$ as a function on $\mathsf{X}$ is a measurable function.
\end{enumerate}
We view it as an continuous linear operator on the space of bounded measurable functions on $B(\mathsf{X})$ \cite{yosida1941}, whose action is denoted by
\begin{equation} \label{markov-op}
    [\mcalP f](x)=\int_{\mathsf{X}} f(y)\mcalP(x,\d y).
\end{equation}
It's evident that $\|\mcalP\|_{B(\mathsf{X})}\le 1$.

The dual space of $B(\mathsf{X})$ is $\mathrm{ba}(\mathsf{X})$, the set of finitely additive set functions \cite{dunford1957}.
The space of signed measures $\mathscr{M}(\mathsf{X})$, which are \textit{countably} additive set functions, and therefore stronger than the requirements of $\mathrm{ba}(\mathsf{X})$, is a closed subspace, and so is the set of probability measures $\mathcal{P}(\mathsf{X})$ (it's a submanifold, to be more precise).
The dual operator $\mcalP^*$ is defined naturally on $\mathrm{ba}(\mathsf{X})$.
\begin{equation} \label{markov-dual}
    \mcalP^* \mu = \mu\mcalP = \int_{\mathsf{X}} \mu(x)\mcalP(x,\cdot).
\end{equation}
However, due to the special form of Markovian kernels, one can show that both $\mathscr{M}(\mathsf{X})$ and $\mathcal{P}(\mathsf{X})$ are closed under the operation of $\mcalP^*$, . Especially, when restricted to $\mathscr{P}$, the action of $\mcalP^*$ returns the probability measure after one transition of the Markov chain.

\section{Numerical scheme}
Recall that $\xi^n$ is the size at birth of the $n$-th along the lineage of observation. 
We notice that in the transition probability density function \eqref{trans-density}, $B(x), g(x)$ always appeared together as their ratio in the Markovian model, so in the following text we will write $S(x)=\frac{B(x)}{g(x)}$. Using $S(x)$ as the primary varialble is equivalent to using $B(x)$, because $g(x)$ is usually known beforehand. 
\begin{equation} \label{trans-tail}
    \BbbP(\xi^{n+1}>y|\xi^n=x)
    =\exp\left( -\int_x^{2y} S(t)\d t \right), \quad \forall y>\frac{x}{2},
\end{equation}
and the Markovian transition kernel determined by this relation is denoted by $\mcalP$.
We will use a finite volume type scheme to implement the GF chain.

\subsection{Finite volume approximation}
First, we shall set up a finite-dimensional approximation of the distributions in $\mathscr{P}(\BbbR_+)$.
Let $a$ be a positive number, denoting a finite range. Let $N_x$ be a positive even integer and $h=N_x^{-1}a$ the mesh size.
Use $x_j=jh, j=0,1,\ldots,N_x$ to denote discrete grid points. The interval $(0,a]$ can then be divided into $N_x$ smaller intervals
\begin{equation}
    E_j\triangleq (x_{j-1}, x_j], \quad j=1,\ldots,N_x.
\end{equation}
They are called \textbf{units}. Also, use $E_{N_x+1}=(a,a+h]$ to denote a unit appended to the right of our grid, and $E_\infty=(a,\infty)$ the infinite tail beyond the grid as a single unit. Therefore, the whole space is partitioned into the disjoint union
$$\BbbR_+=\sum_{j=1}^{N_x} E_j + E_\infty.$$
\begin{figure}[H]
    \centering
    \includegraphics[width=0.5\textwidth]{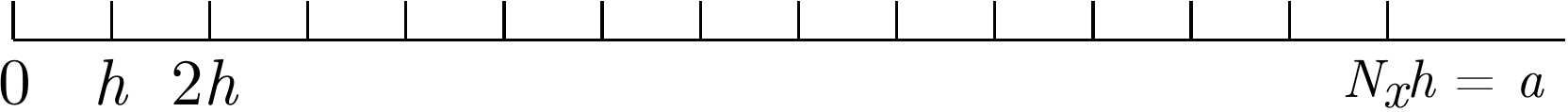}
    \caption{Grid points and units}
\end{figure}

For any probability measure $\mu$ on $\BbbR_+$, we define the following finite volume approximation.
\begin{definition}
    Let $I_{a,h}: \mathscr{P}(\BbbR_+)\to \mathscr{P}(\BbbR_+)$ be the operator on the set of probability distriubtions on $\BbbR_+$, which is defined by
    \begin{equation} \label{fv}
        I_{a,h} \mu = \sum_{j=1}^{N_x} \mu(E_j) \mathcal{U}(E_j)+ \mu(E_\infty)\mathcal{U}(E_{N_x+1}),
    \end{equation}
    where the notation $\mathcal{U}(E)$ refers to the probability measure of a uniform distribution on the set $E$.
    We call it the \textbf{finite volume approximation} of the probability measure $\mu$ with \textbf{range} $a$ and \textbf{mesh size} $h$. 
    We also say a random variable $\zeta$ is a finite volume approximation of another $\eta$ if their distribution laws satisfy 
    $$\mathcal{L}(\zeta)= I_{a,h} \mathcal{L}(\eta).$$
\end{definition}
In the formula \eqref{fv}, $\mathcal{U}(E)$ stands for a measure corresponding to the uniform distribution on the measurable set $E$.
It is evident that $I_{a,h}: \mathscr{P}(\BbbR_+)\to \mathscr{P}(\BbbR_+)$ is a bounded linear projection operator (i.e. $I_{a,h}^2=I_{a,h}$), and that the range of $I_{a,h}$ is finite-dimensional, comprising the space of piecewise uniform measures on $(0,a+h]$.

The operator $I_{a,h}$ is composed of two effects --- granularization and truncation.
\begin{enumerate}
    \item Granularization: the first summation term in \eqref{fv} maps the probability (or ``mass'') of $\mu$ on the unit $E_j$ to a uniform measure over this interval with the same mass, $\forall 1\le j\le N_x$;
    \item Truncation: the second term in \eqref{fv} sends the probability on the infinite tail $E_\infty$ to the unit $E_{N_x+1}$ appended on the right of the grid.
\end{enumerate}
Also, $I_{a,h}$ preserves the total mass on each unit $E_j$ within the grid ($1\le j\le N_x$) as well as the total measure on $E_\infty$, the last of which was sent entirely to $E_{N_x+1}$. 
\begin{gather*}
    [I_{a,h}\mu](E_j)=\mu(E_j), \quad j=1,\ldots,N_x, \\ 
    [I_{a,h}\mu] (E_{N_x+1})=\mu(E_\infty).
\end{gather*}
The range $\ran(I_{a,h})$ is spanned by the $N_x+1$ uniform distributions probabilities
$$\mathcal{U}(E_j),\quad j=1,\ldots,N_x+1,$$
so we can represent such an image with a vector $\bm\mu=(\mu_1,\ldots,\mu_{N_x+1})\in\BbbR^{N_x+1}$ whose components sum to one.
\begin{equation} \label{fv-range}
    \mu_j = [I_{a,h}\mu](E_j), \quad j=1,\ldots,N_x+1.
\end{equation}

With this finite projection operator, we will be able to approximate the GF transition rule \eqref{trans} with a finite matrix, so that computation of an invariant measure is feasible. Also, we will be proving the convergence of this finite approximation as the main part of this paper.

The reason why $N_x$ should be even, and why the probability of the infinite tail $(a,\infty)$ is mapped outside of the grid are purely technical. See the following subsection.

\subsection{Finite transition scheme}
We will define an approximate Markovian transition kernel on $\BbbR_+$, denoted by $\mcalP_{a,h}$, with $a,h$ given in \eqref{fv}.
The range of this probability transition kernel (as a dual operator) will stay in $\ran(I_{a/2,h})$, so that it can be represented by a probability transition matrix, and that the implementation on computers is feasible.

The scheme is devised as follows. 
To distinguish between the two chains, the state variables of this new chain will be denoted by $\{\tilde\xi^n\}$. Suppose we are currently at the state $\tilde\xi^n\in\BbbR_+$.
\begin{enumerate}
    \item Round to the grid: If $\tilde\xi^n\in E_i$, round it to $x_i$, which is right above $\tilde\xi^n$; if $\tilde\xi^n>a$, however, choose $x_i=a$.
    We are doing this because the transition kernel $\mcalP(x_i,\cdot)$ starting at a grid point is easier to approximate than $\mcalP(x,\cdot)$ through quadrature, as will be clarified in the next subsection.
    \item Generate size at fragmentation:
    Define the random variable $\eta$ by
    \begin{equation}
        \BbbP(\eta>y)=\e^{-\int_{x_i}^{y} S(t)\d t}, \quad y\ge x_i.
    \end{equation}
    It models the size right before the fragmentation of a cell with starting mass $x_i$, that is, the distribution of $2\xi^{n+1} | \xi^n=x_i$. We do not sample from this continuous distribution explicitly because it's genuinely hard. Instead, we define $\tilde\eta$ as the finite volume approximation of $\eta$, whose finite dimensional nature allows us to sample directly from it.
    \item Divide and project again: Finally, we define the next state variable $\tilde\xi^{n+1}$ to be the finite approximation of $\tilde\eta/2$. 
\end{enumerate}

Let us discuss several key properties of the transition relation.

\begin{proposition} \label{stay-in-half-grid}
    After one step of transition, the random variable $\tilde\xi^{n+1}$ will be supported on the interval $(0,a/2+h]$.
\end{proposition}
\begin{proof}
    Because $\tilde\xi^{n+1}$ is the finite volume approximation to $\tilde\eta/2,$ which is supported on $(0,(a+h)/2]$, and that $I_{a,h}$ maps the mass on each unit to itself, $\tilde\xi^{n+1}$ must be supported on $(0,a/2+h]$.

    Because we have assumed $N_x$ is even, there is a grid point at $x=a/2$. Thus, the mass of $\tilde\eta/2$ on the rightmost interval $(a/2,(a+h)/2]$ is mapped to a locally uniform distribution on $(a/2,a/2+h]$ under $I_{a,h}$.
\end{proof}

\begin{proposition} \label{fv-on-half-grid}
    $\tilde\xi^{n+1}$ is a finite approximation to $\eta/2$ with range $\frac{a}{2}$ and mesh size $h$.
\end{proposition}
\begin{figure}[th]
    \centering
    \includegraphics[width=0.4\textwidth]{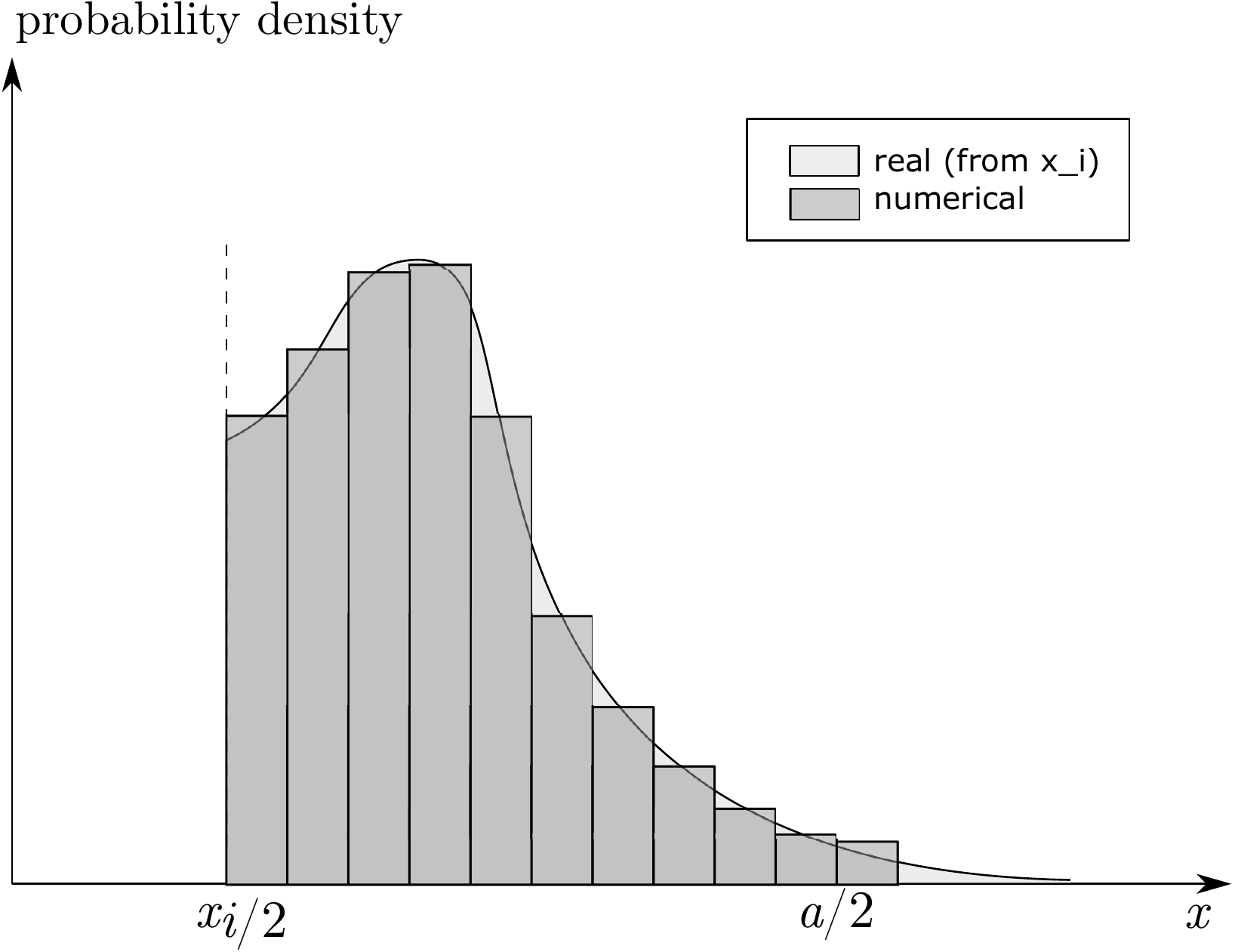}
    \caption{Illustration of Propostion \ref{fv-on-half-grid}. The density function of $\tilde \xi^{n+1}|\tilde\xi^n$ compared with that of $\xi^{n+1}|\xi^n$, both starting at $x_i$.}
\end{figure}
\begin{proof}
    This is because $\tilde\xi^{n+1}$ follows a piecewise uniform distribution over $(0, a/2+h]$, and that the probabilities of $\tilde\xi^{n+1}$ and $\eta/2$ belonging to the corresponding units $E_j$ ($j=1,\ldots, \frac{N_x}{2}+1$) are identical.

    Because $\mathcal{L}(\eta/2)$ is precisely the transitional distribution of the original Markov chain $\xi^{n+1}|\xi^n=x_i$, this assertion also shows that
    $$\tilde \xi^{n+1} \sim I_{a/2,h} \mcalP(x_i,\cdot)$$
    in terms of distribution.
\end{proof}

\begin{proposition} \label{unit-insep}
    The distribution of $\tilde\xi^{n+1}$ is identical for all $\tilde\xi^n\in E_i$. In other words, the units are inseparable under the kernel $\mcalP_{a,h}$.
\end{proposition}
\begin{proof}
    The distribution of $\eta, \tilde\eta$ and therefore $\tilde\xi^{n+1}$ are completely determined by the rounded value $x_i$ from the first step, which is identical for all $\tilde\xi^n \in E_i$.
\end{proof}

As for its transition kernel, denoted by $\mcalP_{a,h}$,
we study the transitional rules of their distribution laws. They are summarized from the scheme into the following lines.
\begin{subequations} \label{num-rule}
\begin{enumerate}
    \item Round $\tilde\xi^{n}$ to the numerical grid, by finding its index $i$.
    \begin{equation}
        i=\max\{ N_x, \lceil \tilde\xi^n/h \rceil \}.
    \end{equation}
    \item Sample $\tilde\eta$, the finite volume approximation to the distribution $2\xi^{n+1}|\xi^n=x_i$.
    \begin{equation}
        \tilde\eta \sim I_{a,h} \mathcal{L}(\eta), \quad
        \frac12 \eta \sim \mcalP(x_i, \cdot).
    \end{equation}
    \item Halve the variable $\eta$, and project again to get $\tilde\xi^{n+1}$.
    \begin{equation}
        \tilde\xi^{n+1} \sim I_{a,h} \mathcal{L}(\tilde\eta/2).
    \end{equation}
\end{enumerate} 
\end{subequations}
Since it's a Markovian transition relation, it must define a kernel $\mcalP_{a,h}$, which, by Proposition \ref{fv-on-half-grid}, satisfies that
$$\ran(\mcalP_{a,h})\subset \ran(I_{a/2,h}).$$
Now we may only consider its behavior on that finite dimensional subspace of probability measures.

Using the finite dimensional representation $\bm\mu\in\BbbR^{N_x/2+1}$ in \eqref{fv-range} to denote any probability measure in $\ran(I_{a/2,h})$, we may express $\mcalP_{a,h}$ in terms of an $(\frac{N_x}{2}+1)\times (\frac{N_x}{2}+1)$ probaility transition matrix $\widetilde P\in\BbbR^{},$ whose elements are defined by the transition probabilities from the $\frac{N_x}{2}+1$ units to themselves,
\begin{equation} \label{ptm}
    \begin{split}
        [\widetilde P]_{ij}\triangleq\tilde p_{ij}&=\BbbP\left(\tilde\xi^{n+1}\in  E_j \mid \tilde\xi^n \in E_i \right) \\ 
        &=\mcalP_{a,h}(x_i, E_j),  \quad i,j=1,2,\ldots,\frac{N_x}{2}+1.
    \end{split}
\end{equation}
We shall refer to this transition kernel as the ``numerical kernel'' henceforth.


\subsection{Computation of the finite kernel}
We define the following quantities first.
Given $\tilde\xi^n\in E_i, i=0,\ldots,\frac{N_x}{2}+1$, we define
\begin{equation}\label{Qk}
    Q_{i,k}=\exp\left( -h\sum_{j=i+1}^{i\vee k} S(x_j) \right)
    \approx\exp\left( -\int_{x_i}^{x_k\vee x_i} S(x)\d x \right), \quad k=0:N_x.
\end{equation}
Clearly, these values are approximations to the complementary distributions of $\eta$ given $\tilde\xi^{n}\in E_i.$ By the relations \eqref{num-rule}, we get
\begin{equation}
    \BbbP(\tilde\xi^{n+1}>x_k|\tilde\xi^n)
    \approx Q_{i,2k}, \quad 0\le k\le \frac{N_x}{2}.
\end{equation}
Taking the difference of these values of $Q_{i,k}$ gives us the transition probabilities over the units $E_k$. Therefore, we have the following approximate expressions for $\tilde p_{ik}$ in \eqref{ptm}.
\begin{equation} \label{kern-disc}
    \tilde p_{ik}=\BbbP(\tilde\xi^{n+1}\in E_k \mid \tilde\xi_n\in E_i )\approx\begin{cases}
        0, & k<i/2, \\
        Q_{i,2k-2}-Q_{i,2k}, & i/2 \le k \le N_x/2, \\
        Q_{i,N_x}, & k=N_x/2+1, \\
        0, & k>N_x/2+1.
    \end{cases}
\end{equation}
\begin{remark}
    Now we explain why only half of the grid points are used. In practice, it may occur we only have finite information about the functions $g(x)$ and $B(x)$, e.g. only at the grid points $x_i$. In this case, it's relatively hard to accurately represent integrals of $S(x)$ on arbitrary intervals. However, in \eqref{Qk}, the upper and lower limits are both grid points, so a right-end approximation to the integral value is quite easy and accurate.
    Because when computing $\tilde{p}_{ik}$, the index of $Q$ reaches $2k$, the values of $Q_{i,k}$ must be available on the whole grid, although the support of $\tilde\xi^n$ is only approximately half of the grid.
\end{remark}
\begin{remark}
    Because $\tilde p_{ik}$ for all values of $i,k$'s can be computed from $Q_{i,k}$, which are further determined by the cumulative sums of $S_j\triangleq S(x_j)$, at most $N_x$ values, their storage may take only $O(N_x)$ space.
\end{remark}

\begin{figure}[th]
    \centering
    \includegraphics[width=0.7\textwidth]{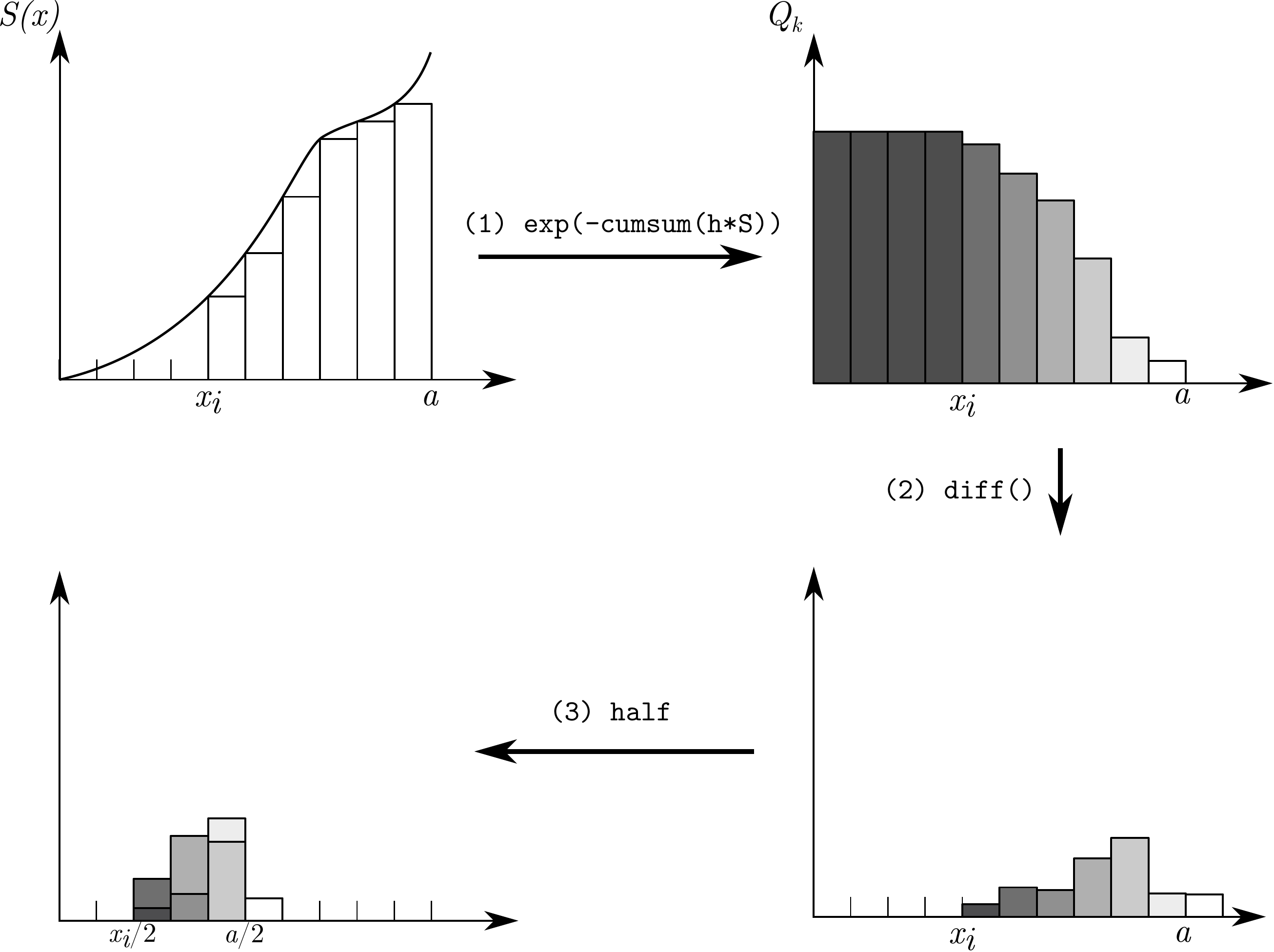}
    \caption{Illustration of the numerical scheme in three steps. Top left: the function values of $S(x)$ drawn as vertical bars; top right: the values $Q_{i,k}$ whose bars are drawn to the left of point $x_k$; bottom right: distribution of $\tilde\eta \approx 2\xi^{n+1}$ by taking the difference of $Q_k$; bottom left: distribution of $\tilde\xi^{n+1}$ through dividing $\tilde\eta$ by 2 and merging adjacent odd and even intervals, the same hue means the same amount of probability.}
\end{figure}

With the values of $\tilde p_{ik}$ obtained, we may use it to generate samples of the Markov chain $\{\tilde\xi^n\}$, or to compute their distributions accurately. Algorithm \ref{gf-alg} generates samples.

\begin{algorithm}[th]
    \caption{Generating samples from the discrete growth-fragmentation chain} \label{gf-alg}
    \KwIn{Previous sample $\tilde\xi^n$; values of $S(x)$ at the gridpoints $\{x_k\}$; grid parameters $a,N_x,h$}
    \KwOut{Next sample $\tilde\xi^{n+1}$}
    $i=\lceil \xi^n/h \rceil$\;
    \If{$i>N_x$}{
        \CommentSty{In this case, $\tilde\eta$ is supported on $(a,\infty)$, so $\tilde\xi^{n+1}$ is supported on only one unit.}\;
        Generate $\tilde\xi^{n+1} \sim \mathcal{U}(a/2,a/2+h]$\; 
    }
    \Else{
        Compute $Q_k$ from \eqref{Qk}\; 
        Compute $\tilde p_{ik}$ from \eqref{kern-disc}\;
        Generate $\tilde\xi^{n+1}$ following the discrete distribution $\BbbP(\tilde\xi^{n+1}\in(x_{k-1},x_k]) = \tilde p_{ik}$\; 
    }
    \Return{$\tilde\xi^{n+1}$}\;
\end{algorithm}

The finite dimensional representation \eqref{fv-range} allows us to represent measures in $\ran(I_{a/2,h})$ with a vector $\bm\mu\in\BbbR^{N_x/2+1}$ and compute the evolution of distribution laws of $\{\tilde\xi^n\}$ as well. The recurrence relation is simply
\begin{equation}
    \bm\mu^{n+1} = \bm\mu^n \widetilde P,
\end{equation}
in which we have computed the transition matrix $\widetilde P$ beforehand.

\section{Convergence theorem}
Here arises a new problem whether this discrete model that we setup provide a correct approximation to the original continuous model, together with the rate of convergence.
In \S\ref{scoot}, we will briefly review a well-known theorem in ergodicity theory. In \S\ref{loo}, we will apply the theory to our continuous Markov chain $\mcalP$ and the finite volume approximation $\mcalP_{a,h}$ to ensure the existence and uniqueness of their invariant measures. as well as apply them to this growth-fragmentation model. Finally, in \S\ref{belle}, we will state our main convergence theorem.

\subsection{Ergodicity} \label{scoot}

Before we introduce the ergodicity theory of Markov chains, we need some basics in functional analysis. For the time being we assume $(\mathsf{X}, \mathfrak{F})$ is a measurable space and that $V:\mathsf{X}\to\BbbR$ is a measurable function such that $V(x)\ge 1,\forall x\in\mathsf{X}$, which will be termed the \textbf{Lyapunov function}.

The \textbf{$V$-bounded measurable functions} are the measurable functions that satisfy
\begin{equation} \label{V-norm}
    \|f\|_V\triangleq \sup_{x\in\mathsf{X}}\frac{|f(x)|}{V(x)}<\infty.
\end{equation}
They form a Banach space denoted by $B_V(\mathsf{X})$ under this norm.
The dual space of $B_V(\mathsf{X})$ is $\mathrm{ba}_V(\mathsf{X})$, the space of finitely additive set functions under which $V$ is integrable \cite{dunford1957} equipped with the \textbf{$V$-total variation norm}
\begin{equation}
    \|\mu\|_{V\textrm{--TV}}\triangleq \<|\mu|,V\>.
\end{equation}
The signed measures and probability measures with finite $V$--TV norms, denoted by $\mathscr{M}_V$ and $\mathscr{P}_V$, are still closed subspaces (submanifolds) of $\mathrm{ba}_V$.

Recall from \eqref{markov-op} and \eqref{markov-dual} that all Markovian kernels can be deemed as linear operators on $B(\mathsf{X})$ and $\mathscr{P}(\mathsf{X})$ by the definition
\begin{align*}
    [\mcalP f](x)&=\int_{\mathsf{X}} f(y)\mcalP(x,\d y). \\
    \mcalP^* \nu &= \nu\mcalP = \int_{\mathsf{X}} \nu(x)\mcalP(x,\cdot).
\end{align*}
With a Lyapunov function $V$ in the model, we need to add an extra assumption
\begin{equation} \label{v-bnd}
    [\mcalP V](x) \le C V(x),
\end{equation}
with which we can show once again that $\mcalP$ and $\mcalP^*$ are still continuous linear operators on $B_V$ and $\mathrm{ba}_V$, under their respective norms.
Additionally, $\mcalP^*$ still keeps $\mathscr{M}_V$ and $\mathscr{P}_V$ invariant because of the Markovian structure, and that it preserves the positivity and the measure of the whole space.

The famous theorem on $V$-uniform ergodicity comes in the following form. It has been stated in such works as \cite{mt1993,mt1994,rosen1995,bax2005,jerison}.
\begin{theorem} \label{erg-thm}
    If the Markovian kernel ${\cal P}$ on a measurable space $(\mathsf{X},\mathfrak{F})$ satisfies the following three conditions,
    \begin{enumerate}
        \item Minorization: There exists a set $C$, a positive constant $\beta$ and a probability measure $\nu$ such that
        \begin{equation} \label{minor}
            {\cal P}(x,E)\ge\beta \nu(E),\quad\forall x\in C, E\in\mathfrak{F}.
        \end{equation} 
        \item Drift condition: There exist constants $\lambda\in(0,1),K>0$ and a measurable Lyapunov function $V:\mathsf{X}\to [1,\infty)$ such that
        \begin{equation} \label{drift}
            {\cal P}V(x)=\int_{\mathsf{X}} V(y)\mcalP(x,\d y)
        \le \begin{cases}
            K, & x\in C,\\
            \lambda V(x), & x\notin C.
        \end{cases}
        \end{equation}
        \item Strong aperiodicity: The probability measure $\nu$ satisfies that $\nu(C)>0$.
    \end{enumerate}
    then it is $V$-uniform ergodic. 

    This means that there exists a probability measure $\pi$ such that for all $V$-bounded measurable functions $f$,
    \begin{equation} \label{V-uni-erg}
        |\mcalP^n f(x)-\<\pi,f\>|\le A \rho^n V(x), \quad\forall x\in \mathsf{X},n\in\mathbb{N}.
    \end{equation}
    $A>0, \rho\in(0,1)$ are computable constants determined by $C,\beta,\nu,\lambda,K$.
\end{theorem}
\begin{remark}
    The drift condition \eqref{drift} ensures the boundedness of $\mcalP$ under the $V$-norm \eqref{v-bnd}: $\|\mcalP\|_V \le \max\{K,\lambda\}$,
    so properties of a general Markovian kernel on the space $B(\mathsf{X})$ also hold on $B_V(\mathsf{X})$.
\end{remark}

There are two direct corollaries from Theorem \ref{erg-thm} in terms of the operator norm of $\mcalP^n$ on the space $B_V(\mathsf{X})$.
\begin{corollary}
    Under the assumptions and notations of Theorem \ref{erg-thm}, the backward operator $\mcalP$ satisfies that
    \begin{equation} \label{pinkie}
        \|\mcalP^n-1\otimes\pi\| \le A\rho^n, \quad  \forall n\in\mathbb{N}.
    \end{equation}
    The LHS stands for the operator norm induced by $\|\cdot\|_V$, and $1\otimes\pi$ is the rank-1 operator $f\mapsto \<\pi,f\>.$
    Equivalently, the dual operator $\mcalP^*$ satisfies that
    \begin{equation} \label{pie}
        \|(\mcalP^*)^n-\pi\otimes 1\|\le A\rho^n, \quad \forall n\in\mathbb{N}.
    \end{equation}
    The LHS stands for the operator norm induced by $\|\cdot\|_{V-\mathrm{TV}}$, and $\pi\otimes 1$ is the rank-1 operator $\mu\mapsto \mu(\mathsf{X})\pi.$
\end{corollary}
\begin{corollary} \label{contract}
    Follow the assumptions and notations in Theorem \ref{erg-thm}.
    Let $\mathscr{M}_0$ denote the space of finite signed measures
    $$\left\{\mu: \<|\mu|,V\><\infty, \mu(\mathsf{X})=0\right\}.$$
    It is an invariant closed subspace of $\mathrm{ba}(\mathsf{X})$ under the operator $\mcalP^*$. 
    Then, the operator $\mcalP^*$ satisfies that
    \begin{equation}
        \|(\mcalP^*|_{\mathscr{M}_0} )^n\|_{V-\mathrm{TV}} \le A\rho^n, \quad \forall n\in\mathbb{N},
    \end{equation}
where $\mcalP^*|_{\mathscr{M}_0}$ stands for the restriction of $\mcalP^*$ to $\mathscr{M}_0$.
\end{corollary}
\begin{remark}
    \eqref{pinkie} is an equivalent statement of \eqref{V-uni-erg}.
    Corollary \ref{contract} shows that the operator $\mcalP^*$ has a \emph{contractive} property when restricted to a closed subspace.

\end{remark}

\subsection{Ergodicity of the GF chain} \label{loo}
Using Theorem \ref{erg-thm}, we can state the ergodicity theorem of the GF chain \eqref{trans} under the assumption of a polynomial growth rate of $S(x)$.
An almost identical has already appeared in \cite{doumic2015}. However, since we have modified their models, we would like to reiterate (and maybe simplify) the proof.

\begin{theorem} \label{gf-erg-thm}
    Assume $S(x)\in C[0,\infty)$, and that there exists constants $m,M,\alpha>0$ and $X_0\ge 0$ such that
    \begin{equation} \label{S-assum}
        m x^{\alpha-1}\le S(x)\le Mx^{\alpha-1}, \quad\forall x\ge X_0.
    \end{equation}
    Then, the Markovian kernel $\mcalP$ \eqref{trans} of the GF chain is $V$-uniform ergodic. Moreover, one may set
    \begin{equation} \label{lyap}
        V(x)=\exp\left( \int_0^x S(t)\d t \right)
    \end{equation}
    as the Lyapunov function in Theorem \eqref{erg-thm}.
\end{theorem}
\begin{proof}
We will verify the three conditions in Theorem \ref{erg-thm}.

Let $C=[0,X]$ for some $X>0$, and $\nu$ the uniform measure on $[X/2,X]$.
The minorization and strong aperiodicity conditions are easy to verify, so we only sketch them without detailed discussion. Using the expression \eqref{trans-tail}, we get
\begin{gather*}
    \mcalP(x, [a,b])= \e^{-\int_{x}^{2a} S(t)\d t}-\e^{-\int_{x}^{2b} S(t)\d t} \ge \e^{-\int_{X}^{2X} S(t)\d t} 2(b-a) \triangleq \beta(b-a), \quad \forall x\in [0,X], \\
    \Rightarrow \mcalP(x, [0,X])\ge \e^{-\int_{X}^{2X} S(t)\d t} X \triangleq \delta>0,
\end{gather*}
which are quite loose bounds.

Next, we verify the drift condition with direct calculation. First we write out the density function \eqref{trans-density} in terms of $V$.
$$p(x,y)=\frac{\pt}{\pt y} \left[\exp\left( -\int_{x}^{2y} S(t)\d t \right)\right]
=\frac{2V(x)V'(2y)}{V^2(2y)}, \quad y>\frac{x}{2}.$$
Substitute it into the formula of $\mcalP V(x)$,
$$\mcalP V(x)=\int_{x/2}^\infty \frac{2V(x)V(y)V'(2y)}{V^2(2y)} \d y 
=V(x) \int_{x/2}^\infty \frac{2V(y)V'(2y)}{V^2(2y)} \d y.$$
Rewrite $V(y)$ with $S(y)$ and extract the $V(x)$ factor from the integral,
$$\mcalP V(x)=V(x) \int_{x/2}^\infty 2S(2y)\e^{-\int_y^{2y} S(t)\d t} \d y.$$
Assume that $x\ge 2X_0$, and use the assumption \eqref{S-assum} to bound $S(x)$. We find that the negative exponent $\int_y^{2y} S(t)\d t$ still grows like a polynomial of $y$.
$$V(x) \int_{x/2}^\infty 2S(2y)\e^{-\int_y^{2y} S(t)\d t} \d y
\le V(x) \int_{x/2}^\infty 2M(2y)^{\alpha-1} \e^{-\int_y^{2y} mt^{\alpha-1}\d t}\d y.$$
And finally, through some basic calculus, we get
\begin{equation} \label{twilight}
    \mcalP V(x) \le C_1 V(x) \e^{-C_2(x/2)^\alpha},
\end{equation}
with the constants defined as
$$C_1=\frac{2^\alpha M}{(2^\alpha-1)m}, \quad C_2=m\alpha^{-1}(2^\alpha-1),$$
both of which fully determined by $m,M,\alpha$ in \eqref{S-assum}.

Clearly, as $x\to\infty$, the coefficient $C_1\e^{-2^{-\alpha}C_2x^\alpha}\to 0$ rapidly, so the second part of \eqref{drift} holds as long as we choose $X$ sufficiently large, e.g. $X=2X_0+1$. On the other hand, it is obvious that $\mcalP V(x)$ is bounded by a uniform constant $K=K(X)>0$ for all $x\in[0,X]$:
$$\mcalP V(x) \le V(X) \int_0^\infty 2S(2y) \e^{-\int_y^{2y} S(t)\d t} \triangleq K, \quad \forall x\in[0,X].$$
Therefore, the conditions are all met, and applying Theorem \ref{erg-thm} gives us the desired property.
\end{proof}

From now on we will denote the unique invariant measure of the GF chain by $\pi$.

In the last step of the proof, the rapid decay of $V^{-1}(y)$ provides us with the rapid decay of $p(x,y)$ in y, and hence the strong bound \eqref{twilight} on $\mcalP V$ for large $x$.
Therefore, we may vaguely believe that the Markov chain determined by $\mcalP$ mostly stays in a finite range. This is a useful result that we shall use later in the proof of convergence. We state it as a lemma below.
\begin{lemma} \label{tail-bnd}
With the GF kernel $\mcalP$ defined by \eqref{trans}, the Lyapunov function $V$ defined by \eqref{lyap}, and under the assumptions \eqref{S-assum}, the following holds for every $x'\in\mathbb{R}$.
\begin{equation}
    \int^\infty_{x'} p(x,y) V(y)\d y \le C_1 V(x) \e^{- C_2 (x')^\alpha},
\end{equation}
where $C_1,C_2>0$ are constants determined by $m,M,\alpha,X_0$.
\end{lemma}
\begin{proof}[Proof]
    The relation has been guaranteed by \eqref{twilight} for $x'\ge X_0$. For $0\le x'< 2X_0$, however, because of the continuity of $S$ and the positivity of $\e^{-C_2 x'^\alpha}$, we can always increase $C_1$ (with a multiplicative constant determined by $X_0$) to make the relation hold uniformly for $x'\in[0,2X_0]$.
\end{proof}

\subsection{Convergence of numerical chain} \label{belle}
Now we are at a point to assert that the numerical kernel $\mcalP_{a,h}$ defined through \eqref{ptm} generates a discrete invariant measure $\pi_{a,h}$ over $\BbbR_+$ that converges to the real invariant GF distribution $\pi$.
The convergence $\pi_{a,h}\to\pi$ will be coupled with the convergence of $\mcalP_{a,h}\to\mcalP$ as is seen in the following statement.
\begin{theorem} \label{main-thm}
    Let the transition kernels ${\cal P}$ and ${\cal P}_{a,h}$ be defined by \eqref{trans} and \eqref{ptm}, respectively. 
    Assume \eqref{S-assum} holds, and that there exists $L>0,\beta>0$ such that
    \begin{equation} \label{smooth-assum}
        \begin{cases}
            S(x)\in C^1[0,\infty), \\ 
            |S'(x)|\le L (1+x^\beta), & \forall x>0.
        \end{cases}
    \end{equation}
    When the parameters satisfy
    \begin{equation} \label{param-constraint}
        \begin{cases}
            h\le 1, \\
            a>2X_0, \\
            h=O(a^{-\alpha+1}), \\
            a\ge 3h,
        \end{cases}
    \end{equation}
    there exist constants $C,k>0$ (determined by $M,m,\alpha,X_0,L,\beta$) such that for all $x\in\mathbb{R}_+$,
    \begin{equation}
        |{\cal P}_{a,h}f(x)-\mcalP_f(x)|\le C V(x) (h+ah+\e^{-ka^\alpha}).
    \end{equation}
    Equivalently, the operator norm of ${\cal P}-{\cal P}_{a,h}$ on the space $B_V(\mathbb{R}_+)$ satisfies
    \begin{equation} \label{main-result}
        \|{\cal P}-{\cal P}_{a,h}\|_V \le \tilde C (h+ah+\e^{-ka^\alpha}).
    \end{equation} 
\end{theorem}
\begin{remark}
    The assumption \eqref{smooth-assum} asserts the smoothness of $S(x)$, which will be used later in the proof. It is natural because the finite volume discretization \eqref{fv}, we are essentially performing numerical quadrature of the transitional density $p(x,y)$ over $y$ (i.e. integrating it over the units $E_j$ and then using them as discrete masses), whose convergence estimate requires the differentiation of $S$.
\end{remark}
\begin{remark}
    The three terms in the error bound to the RHS of \eqref{main-result} can be explained. The first two are of order $O(h)$, which is reasonable since our finite volume discretization of moving a point $x$ to its neighbour $x_i$ by a distance no more than $h$ is generally first order. The last one is rapidly decaying in $a$, and can be ignored.
    In practice, we may first choose a moderately large $a$ enough to extinguish the third term, and then let $h\to 0$ to deal with the first two terms.
\end{remark}

\section{Proof of the main results}
This section is dedicated to the complete proof of Theorem \ref{main-thm}. In \S\ref{celestia}, we will discuss the overall approach. In \S\ref{bloom} and \S\ref{sweet} we will estimate two terms in the error bound respectively. In \S\ref{apple}, we will present the final upper bound.

\subsection{Triangle inequality argument} \label{celestia}
The invariant measures $\pi_{a,h}$ and $\pi$ can be seen as eigenvectors of the dual operators $\mcalP_{a,h}^*$ and $\mcalP^*$ associated with the eigenvalue 1. When the Markov kernel is $V$-uniform ergodic, the eigenvalue is simple, meaning the kernel is contractive on the quotient space over this distribution, as seen from Corollary \ref{contract}.
The convergence $\pi_{a,h}\to\pi$ can be proved with $\mcalP_{a,h}\to\mcalP$ from this perspective.

To make this approach work, we should prove the following two facts:
\begin{enumerate}
    \item Stability: an invariant measure $\pi_{a,h}$ exists for $\mcalP_{a,h}$, and is uniformly bounded for all $a,h$ satisfying \eqref{param-constraint};
    \item Consistency: as $a\to\infty, h\to 0$ at some specific rate, the numerical kernel $\mcalP_{a,h}$ converges to the real kernel $\mcalP$ in operator norm;
\end{enumerate}
These two properties will ensure $\mcalP_{a,h}\to\mcalP \Rightarrow \pi_{a,h}\to\pi$.

First, we have the eigenvalue relations
\begin{equation}
    \begin{cases}
        \pi_{a,h}\mcalP_{a,h}=\pi_{a,h}, \\
        \pi \mcalP =\pi.
    \end{cases}
\end{equation}
Subtracting them, using $\pi_{a,h}\mcalP$ as an intermediate amount and applying the triangle inequality, we get
$$\|\pi_{a,h}-\pi\|
=\|\pi_{a,h}\mcalP_{a,h}-\pi\mcalP\| 
\le \|\pi_{a,h}\mcalP_{a,h}-\pi_{a,h}\mcalP\|
+\|\pi_{a,h}\mcalP-\pi\mcalP\|.$$
Because $\pi-\pi_{a,h}$ is in $\mathscr{M}_0$ (see Corollary \ref{contract} for a definition), we have
\begin{equation}
    (1-\|\mcalP^*|_{\mathscr{M}_0}\|) \|\pi_{a,h}-\pi\| \le \|\mcalP_{a,h}^*-\mcalP^*\| \cdot \|\pi_{a,h}\|.
\end{equation}
Corollary \ref{contract} makes it feasible to assume that $\|\mcalP^*|_{\mathscr{M}_0}\|<1$, because if not we may as well replace $\mcalP^*$ with $(\mcalP^*)^n$, where $n$ is sufficiently large.
The coefficient on the LHS is thus positive, and we can divide it to the RHS.
\begin{equation} \label{luna}
    \|\pi_{a,h}-\pi\| \le \frac{\|\pi_{a,h}\|}{1-\|\mcalP^*|_{\mathscr{M}_0}\|}\|\mcalP_{a,h}^*-\mcalP^*\|.
\end{equation}
Because the RHS converges to 0, so does the RHS.
Therefore, when compatibility and stability hold, $\mcalP_{a,h}\to\mcalP$ implies $\pi_{a,h}\to\pi$.

The consistency condition is the central subject of this section.
The stability condition can be reduced to consistency in our case by converting it to the $V$-uniform ergodic estimate of the numerical kernel $\mcalP_{a,h}$. With Theorem \ref{erg-thm} and its corollaries, in order to get a bound for $\pi_{a,h}$, we will only need uniform constants in the sufficient conditions therein. These conditions can be verified by the convergence of $\mcalP_{a,h}\to\mcalP$.
For any measurable set $E$, denoting by $I_E(x)$ the indicator function of $E$,
$$I_E(x)=\begin{cases}
    1, & x\in E, \\ 0, & x\notin E,
\end{cases}$$
we get
$$\mcalP_{a,h}(x,E)=\mcalP_{a,h} I_E(x)
\ge \mcalP I_E(x)- \|\mcalP-\mcalP_{a,h}\|_V V(x).$$
Choose the same small set $C=[0,X]$ and $\nu=\mathcal{U}[X/2,X]$. If $\mcalP_{a,h}\to\mcalP$, we can easily find uniform constants $\beta',\delta'$ in the minorization and strong aperiodicity conditions.
Moreover,
$$\mcalP_{a,h} V(x) \le \mcalP V(x)+ \|\mcalP-\mcalP_{a,h}\|_V V(x),$$
so by estimating $\|\mcalP-\mcalP_{a,h}\|_V$, we can obtain uniform bounds $\lambda',K'$ in the drift condition.

As long as the $V$-uniform ergodicity of $\mcalP_{a,h}$ is proved, the uniform boundedness of $\|\pi_{a,h}\|_V$ will follow easily from \eqref{V-uni-erg} by setting $n=0$,
\begin{align*}
    &\sup_x \frac{|f(x)-\<\pi_{a,h},f\>|}{V(x)} \le A(\beta',\delta',\lambda',K')  \\
&\Rightarrow 
\|\pi_{a,h} -I\|_V \le A \\
&\Rightarrow \|\pi_{a,h}\|_V \le A+1,
\end{align*}
and it will be clear from \eqref{luna} that convergence $\pi_{a,h}\to \pi$ is obtainable from $\mcalP_{a,h}\to\mcalP$.

\begin{remark}
The triangle inequality technique we use above is directly found in \cite{jin_ergodicity_2022, ye_error_2022}.
It is identical to that in Remark 6.3 of \cite{mattingly2010}.
Similar methods were also used in \cite[Thm. 3.1]{shardlow2000} and \cite[Thm. 7.3]{mattingly2002} for an estimation of the distance between two invariant measures through the distance between their respective processes.
\end{remark}

We are going to compute the error in operator norm $\|\mcalP-\mcalP_{a,h}\|_V$ through its definition \eqref{markov-op}, that is, by estimating
\begin{equation} \label{error}
    \sup_{x\in\BbbR_+} \frac{|\mcalP f(x) - \mcalP_{a,h} f(x)|}{V(x)} 
    =\sup_{x\in\BbbR_+} \frac{1}{V(x)} \left| \<\mcalP(x,\cdot) ,f\>- \<\mcalP_{a,h}(x,\cdot), f(x)\> \right|.
\end{equation}
for all $f\in B_V, \|f\|_V\le 1$.

According to Propositions \ref{stay-in-half-grid}--\ref{unit-insep}, the transitional distribution $\mcalP_{a,h}(x,\cdot)$ is identical to $\mcalP_{a,h}(x_i,\cdot)$, the latter of which is closely related to the continuous distribution $\mcalP(x_i,\cdot)$ by the finite volume approximation.
Therefore, to compute \eqref{error}, we can use the distribution $\mcalP(x_i,\cdot)$ as a medium, which splits it into two parts.
\begin{equation} \label{error-part}
    |\mcalP f(x)-\mcalP f(x_i)|
+ |\mcalP f(x_i)-\mcalP_{a,h} f(x)|.
\end{equation}
The first part is a change of the continuous kernel in the first argument $x$ (i.e. the ``starting point'' of the transition), and the second is the discretization of a fixed transition probability from $x_i$. Let's estimate them in the following two subsections respectively.

\begin{remark}
    More strictly speaking, our numerical scheme only \textit{approximates} the discrete transition probabilities $\tilde p_{ij}$ by the values of $Q_{i,k}$ (recal \eqref{Qk}). However, for simplicity we will now assume that the discrete transition probabilities $\tilde p_{ij}$ are computed accurately, i.e. without quadrature error. This is negligible because as $h\to 0$, their values (and hence the difference of their corresponding transition kernels) tend together at a rate of $O(h)$, which fits in the RHS error estimate of \eqref{main-result}.
\end{remark}

\subsection{Estimating a change in the starting point} \label{bloom}
The following lemma will provide the estimate of the first summand of \eqref{error-part}.
\begin{lemma} \label{start-err}
    Suppose that $S(x)$ is continuous on $[0,\infty)$ and that the assumption \eqref{S-assum} holds for $x>X_0$. Then there exists a uniform constant $C_3$ that is only determined by $M,m,\alpha,X_0$ such that for sufficently small $h$ (i.e. $h\le 1$) and any $x\in\mathbb{R}_+$,
    \begin{equation}
        \|{\cal P}(x, \cdot)-{\cal P}(x+h,\cdot)\|_{V-\mathrm{TV}}\le C_3 h V(x+h).
    \end{equation}
\end{lemma}
\begin{proof}
We can bound the LHS by
$$\int_0^\infty |p(x,y)-p(x+h,y)| V(y)\d y.$$
The function $p(x,\cdot)$ is supported on $[x/2,\infty)$, while $p(x+h,\cdot)$ on $[(x+h)/2,\infty)$.
Therefore we split the domain into two parts:
\begin{equation}
    (\text{I})=\int_{x/2}^{(x+h)/2} |p(x,y)-p(x+h,y)| V(y)\d y
=\int_{x/2}^{(x+h)/2} p(x,y) V(y)\d y,
\end{equation}
because $p(x+h,y)$ is zero for $y<\frac{x+h}{2}$, and
\begin{equation}
    (\text{II})=\int_{(x+h)/2}^\infty \int_{x/2}^{(x+h)/2} |p(x,y)-p(x+h,y)| V(y)\d y,
\end{equation}

For the first part, we have
$$(\text{I})
=\int_{x/2}^{(x+h)/2} p(x,y) V(y) \d y  \\
=\frac{h}{2} p(x,s)V(s),$$
with $2s\in[x,x+h]$ by the Mean Value Theorem. Express all of the quantities in terms of the function $S(x)$.
$$(\text{I})=h \cdot S(2s) \cdot \e^{-\int_x^{2s} S(t)\d t} \e^{\int_0^s S(t)\d t}
=h V(x) S(2s)\e^{-\int_s^{2s} S(t)\d t}.$$
During the proof of Theorem \ref{gf-erg-thm}, we obtained the following bound.
$$\int_{s}^{2s} S(t)\d t \ge C_2 s^{\alpha}.$$
See the proof there for the explicit expression of $C_2$.
Using it once again we get
\begin{equation} \label{flutter}
(\text{I})\le h V(x) \e^{-2^{-\alpha}C_2 x^\alpha} \sup_{[x,x+h]}S .
\end{equation}

For the second part, we notice the following relation from the form of $p(x,y)$ \eqref{trans-density}.
$$p(x+h,y)=p(x,y)\cdot \e^{\int_x^{x+h} S}, \quad y>\frac{x+h}{2}.$$
Therefore,
\begin{align*}
    (\text{II})
    &=\int_{(x+h)/2}^\infty p(x+h,y) V(y) \left|1-\e^{-\int_x^{x+h} S(t)\d t}\right|\d y \\
    &\le \int_x^{x+h} S(t)\d t
\int_{(x+h)/2}^\infty p(x+h,y) V(y)\d y
\end{align*}
We have used the basic inequality $1-\e^{-x}\le x,\forall x\ge 0.$
Then, we bound the first integral by a supremem and apply Lemma \ref{tail-bnd} to the second integral to get
\begin{equation} \label{shy}
    (\text{II}) \le C_1 h V(x+h) \e^{-2^{-\alpha}C_2 (x+h)^\alpha}\sup_{[x,x+h]} S.
\end{equation}

Since $S(x)$ grows like a polynomial, it is obvious that
\begin{equation} \label{rarity}
    \sup_{x\ge 0} \left[ \e^{-kx^\alpha}\sup_{[x,x+h]}S \right]<\infty, \quad \forall \alpha>0, k>0.
\end{equation}
Suppose further that $h\le 1$, and we know that \eqref{rarity} is bounded by a constant.
Also, $V(x)\le V(x+h)$ because $S\ge 0$, so the sum of \eqref{flutter} and \eqref{shy} can be further bounded by 
$$(\text{I})+(\text{II})\le C_3 h V(x+h),$$
where $C_3$ is a constant depending on $M, C_1$ and $C_2$ only, and hence determined by $M,m,\alpha,X_0$. The proof is then complete.
\end{proof}

\subsection{Estimating the numerical error} \label{sweet}
Next, we estimate the second term in \eqref{error-part}. Because $\mcalP_{a,h}(x,\cdot)=\mcalP_{a,h}(x_i,\cdot),$ it is equivalent to estimating
\begin{equation} \label{applejack}
    \mcalP f(x_i)-\mcalP_{a,h} f(x_i)=\int_{x_i/2}^\infty p(x,y)f(y)\d y-\int_{\mathbb{R}} {\cal P}_{a,h}(x_i,\d y) f(y)
\end{equation}
for each $x_i, \|f\|_V\le 1$.

Now the extra smoothness assumption \eqref{smooth-assum} in the main theorem comes into play. Recall we asserted that there exist constants $L,\beta\ge 0$ such that 
$$\begin{cases}
        S(x)\in C^1[0,\infty), \\ 
        |S'(x)|\le L (1+x^\beta), & \forall x>0.
\end{cases}$$
Using the smoothness assumption, the following lemma is obtained to describe the smoothness of $p(x,y)$ over $y$, which can be understood as some type of ``Lipschitz'' estimate.
\begin{lemma} \label{p-lip}
Under the assumptions \eqref{S-assum} and \eqref{smooth-assum}, the probability density satisfies that
\begin{equation}
    |p(x,y)-p(x,y+h)|\le C_4 h V(x)\frac{(1+y^\gamma)}{V(2y)} ,        
\end{equation}
for all $x>0, h\ge 0$ and $y>\frac{x}{2}$. $C_4=C_4(M,m,\alpha,X_0)$ is a constant independent of $x,y,h$, and $\gamma=\beta\vee 2(\alpha-1).$
\end{lemma}

\begin{proof}
Differentiate $p(x,y)$ with respect to $y$.
$$\frac{\pt p(x,y)}{\pt y}=4[S'(2y)-S^2(2y)]\e^{-\int_x^{2y} S(t)\d t}.$$
So with the fundamental theorem of calculus, applying assumptions \eqref{S-assum} and \eqref{smooth-assum}, we get
\begin{align*}
    |p(x,y+h)-p(x,y)|
&\le \int_{y}^{y+h} 4|S'(2z)-S^2(2z)|\e^{-\int_x^{2z}S(t)\d t} \d z \\
&\le \frac{V(x)}{V(2y)} \int_{y}^{y+h} 4[L+L(2z)^{\beta}+M(2z)^{2(\alpha-1)}+C(X_0)] \d z.
\end{align*}
Evidently The latter term is of polynomial order $O( h(1+y^\gamma))$, so we can use $1+y^\gamma$ to bound it, above the constant $C_4$ is only dependent on $L, M,\alpha,X_0,$ which is independent from $x,y,h$.
\end{proof}

The next lemma provides the estimates of \eqref{applejack}. The logic behind this result and its assumptions will be clarified in the proof.
\begin{lemma} \label{disc-err}
    Assume that \eqref{S-assum} and \eqref{smooth-assum} holds. Assume further that $h\le 1$ and $3h\le a$. Then, the numerical error satisfies that
    \begin{equation}
        |\mcalP f(x_i)- \mcalP_{a,h} f(x_i)|
        \le C_8 V(x_i) [ah+\e^{-C_6 a^\alpha}],
    \end{equation}
    where the constants $C_8,C_6$ depend only on $M,m,\alpha,X_0,L,\beta$.
\end{lemma}
\begin{proof}
Because $\mcalP_{a,h}(x_i,\cdot)$ is a combination of piecewise constant probability measures, we can write out $\mcalP_{a,h} f(x_i)$ explicitly with the rules \eqref{num-rule} and \eqref{ptm}.
$${\cal P}_{a,h} f(x_i)=
\sum_{j=\lceil i/2\rceil}^{N_x/2} {\cal P}(x_i, E_j) \intbar_{x_{j-1}}^{x_j} f
+{\cal P}(x_i, (a/2,\infty)) \intbar_{a/2}^{a/2+h} f.
$$
Here, $\smallintbar_A f$ stands for the integral average of a function $f$ over the set $A$ and $\{E_j\}$ are the finite volumes units in the grid.
We are using the exact transition probabilities to compute $\tilde p_{ij}$ because of Proposition \ref{fv-on-half-grid}.

We partition the integral region in $\mcalP f(x_i)$ into the same intervals $(x_{j-1},x_j]$ for $j=1:N_x/2$ and $(a/2,\infty)$.
Then, the error \eqref{applejack} may be split into two parts again. The first is
\begin{equation}
    (\text{III})=
\sum_{j=1}^{N_x/2} \left[
    \int_{x_{j-1}}^{x_j} f(y) p(x_i,y)\d y - \intbar_{x_{j-1}}^{x_j} f(y)\d y \int_{x_{j-1}}^{x_j} p(x_i,y)\d y
\right],
\end{equation}
which stands for the discretization error on the intervals over $[0,a/2]$, and is believed to be $O(h)$.
The second is
\begin{equation}
    (\text{IV})=\int_{a/2}^\infty p(x,y)f(y)\d y
-\intbar_{a/2}^{a/2+h} f(y)\d y\int_{a/2}^\infty p(x,y)\d y .
\end{equation}
and the truncation error of $[a/2,\infty),$ which decays exponentially in $a$. 

\paragraph{Discretization error}

If $x_j>x_i/2$, $p(x_i,y)f(y)$ is integrated on the whole interval, and then we have the following termwise error of (III).
\begin{align*}
    &\int_{x_{j-1}}^{x_j} f(y) p(x_i,y)\d y - \intbar_{x_{j-1}}^{x_j} f(y)\d y \int_{x_{j-1}}^{x_j} p(x_i,y)\d y \\ 
    =\;&\int_{x_{j-1}}^{x_j} f(y) p(x_i,y)\d y - \int_{x_{j-1}}^{x_j}\frac1h f(y)\d y \int_{x_{j-1}}^{x_j} p(x_i,y)\d y \\ 
    =\;& \int_{x_{j-1}}^{x_j} f(y) \left(p(x_i,y)-\intbar_{x_{j-1}}^{x_j} p(x_i,z)\d z\right)\d y \triangleq E_j.
\end{align*}
If $x_{j-1}\ge x_i/2$, $p(x_i,y)$ is smooth over the whole interval $(x_{j-1},x_j]$. By Lemma \ref{p-lip} and the Mean Value Theorem, for all $y\in (x_{j-1},x_j],$
\begin{align*}
    |E_j|&=\left|p(x_i,y)-\intbar_{x_{j-1}}^{x_j} p(x_i,z)\d z\right| \\
    &=|p(x_i,y)-p(x_i,\eta)|
    \le C_4 \frac{V(x_i)}{V(2x_{j-1})} (1+x_j^\gamma) |y-\eta|.
\end{align*}
In the last formula, we increased the bound from Lemma \ref{p-lip} by shrinking $V(2y)$ in the denominator to $V(2x_j-1)$ and expanding $1+y^\gamma$ to $1+x_j^\gamma$.
Moreover, we have $|f(y)|\le V(y)\le V(x_j)$.
So an upper bound for the termwise error $E_j$ is
\begin{equation} \label{ej}
    |E_j|\le C_4 V(x_i) \frac{(1+x_j^\gamma) V(x_j)}{V(2x_{j-1})} 
\int_{x_{j-1}}^{x_j} |y-\eta| \d y.
\end{equation}
Assume $h\le 1$, and it's evident that the factor $\frac{(1+x_j^\gamma) V(x_j)}{V(2x_{j-1})} $ is uniformly bounded for all $j$ and $h\le 1$, again because of the rapid decay of an exponential term:
$$\sup_{x\in\BbbR, h\le 1} \frac{[1+(x+h)^\gamma]V(x+h)}{V(2x)}=C(M,m,\alpha,X_0,\gamma)<\infty.$$
Therefore, we can bound $E_j$ further with
$$|E_j|\le C_5 V(x_i) h^2.$$

If $i$ is even, every interval in the sum (III) either satisfies $x_j\ge x_i/2$, or is completely zero for both kernels. Therefore, we may sum up the $E_j$'s from $j=i/2$ to $j=N_x/2$ to get
$$(\text{III})\le \frac12 C_5 V(x_i) ah.$$

If $i$ is odd, there is an interval $[x_{\lfloor i/2\rfloor}, x_{\lceil i/2 \rceil}]$ on which $p(x_i,\cdot)$ is not smooth.
But this term is uniformly of order $h$, so it doesn't really affect the bound.
Let $i=2j-1\Leftrightarrow j=\lceil i/2 \rceil$, and the term $E_j$ now equals
$$\int_{x_i/2}^{x_j} f(y)p(x_i,y)\d y 
    -\int_{x_{j-1}}^{x_j} f(y)\d y \intbar_{x_{j-1}}^{x_j} p(x_i,z)\d z.$$
The first term has already been estimated in the proof of Lemma \ref{start-err} as the quantity (I). By replacing the variable $x$ there with $x_i/2$, we immediately get
$$\left| \int_{x_i/2}^{x_j} f(y)p(x_i,y)\d y  \right|\le \int_{x_i/2}^{x_j} p(x_i,y) V(y)\d y
\le (\tilde C-1) hV(x_i)$$
The second term is directly estimated by its upper bound. 
$$\left|\int_{x_{j-1}}^{x_j} f(y)\d y \intbar_{x_{j-1}}^{x_j} p(x_i,z)\d z \right|
< h\sup_{[x_{j-1},x_j]}|f| \cdot 1 <hV(x_j)\le hV(x_i).$$
So we have the estimate $|E_j|\le \tilde Ch V(x_i)$
Summing all the rest of the terms, we get
$$(\text{III})\le \tilde C_5 V(x_i) ah.$$

As a summary, no matter what the parity of $i$ is, we all have an upper bound of $O(ah)$ for (III) if $h\le 1$.
\begin{equation} \label{rainbow}
    (\text{III})\le C_5 V(x_i) ah,
\end{equation}
where $C_5$ only depends on $M,m,\alpha,X_0,L,\beta$.

For the term (IV), we only have to estimate the infinite tails.
Apply Lemma \ref{tail-bnd} to the first term.
$$\int_{a/2}^\infty p(x_i,y) f(y)\d y 
\le C_1 V(x_i)\e^{-C_2 2^{-\alpha}a^\alpha}.$$
Transform the second term a little by writing out $p(x_i,y)$.
\begin{align*}
    \left| \frac{1}{h}\int_{a/2}^{a/2+h} f(y)\d y\int_{a/2}^\infty p(x,y)\d y \right|
    &\le\frac{1}{h} \int_{a/2}^{a/2+h} V(y)\d y
    \int_{a/2}^\infty 2S(2y)\frac{V(x_i)}{V(2y)}\d y \\
    &=V(x_i)\cdot \intbar_{a/2}^{a/2+h} \e^{-\int_y^a S(t)\d t}\d y
    \int_{a/2}^\infty 2S(2y) \e^{-\int_a^{2y} S(t)\d t} \d y.
\end{align*}
The third factor is always 1 because it's just the kernel $p(a,\cdot)$. The second factor can be bounded by
$$\intbar_{a/2}^{a/2+h} \e^{-\int_y^a S(t)\d t}\d y
\le \exp\left( \int_{a/2+h}^a S(t)\d t \right)
\le \exp\left[-m\alpha^{-1}(a^\alpha-(a/2+h)^\alpha)\right].$$
It decays with exponential order in $a$, as long as $h$ is not too large. If we assume that $3h\le a$, we may bound it by some $\e^{- \tilde C_6 a^\alpha}.$

Summing two terms up and taking the maximum, we can get the following bound for (IV).
\begin{equation} \label{dash}
    (\text{IV})\le
    C_7 V(x_i) \e^{-C_6 a^\alpha},
\end{equation}
with the constants $C_6,C_7$ only depending on $M,m,\alpha,X_0$.

Summarizing \eqref{rainbow} and \eqref{dash}, we get that for all $\|f\|_V\le 1$,
$$|\mcalP f(x_i)- \mcalP_{a,h} f(x_i)|
\le (\text{III})+(\text{IV})
\le C_8 V(x_i) [ah+\e^{-C_6 a^\alpha}]$$
as desired.
\end{proof}

\subsection{The overall estimate} \label{apple}
Now we are ready to provide the ultimate estimate of $\mcalP-\mcalP_{a,h}.$

First, let's deal with trivial case where $x>a$. The discrete transition probability $\mcalP_{a,h}(x,\cdot)$ is supported on $(\frac{a}{2},\frac{a}{2}+h]$ by Proposition \ref{fv-on-half-grid}.
We may estimate the \textit{sum} of the two terms in the error \eqref{error} directly through Lemma \ref{tail-bnd}.
\begin{gather*}
    {\cal P}_{a,h}f(x)=\intbar_{a/2}^{a/2+h} f 
\le V(a/2+h)\le V(x) \e^{-\int_{a/2+h}^x S}, \\
{\cal P}f(x)\le C_1 V(x) \e^{-C_2 (x/2)^\alpha}.
\end{gather*}
When $h\le 1$ and $a\ge 3h$, Both are bounded by $V(x)\e^{-kx^\alpha}.$ (This $k$ may be different from the previous $C_6$. We choose the smaller one.) Therefore,
\begin{equation}
    |{\cal P}_{a,h}f(x)-{\cal P} f(x)|\le C V(x) \e^{-k x^\alpha},\quad x>a,
\end{equation}
and it can certainly be controlled by $\e^{-k a^\alpha}.$

Suppose that $x\in(x_{i-1},x_i].$ Adding the results of the Lemmas \ref{start-err} and \ref{disc-err} together, we get
\begin{equation} \label{sparkle}
    |{\cal P}f (x_i)-{\cal P}_{a,h}f(x)| \le 
    C V(x_i) [h+ah+ \e^{-k a^\alpha}].
\end{equation}
where the constants $C,K$ does not depend on $x, h$ as long as $h\le 1, 3h\le a$. We want the ratio of the RHS of \eqref{sparkle} and $V(x)$ to be bounded, but $V(x_i)$ is still larger in magnitude than $V(x)$.
More precisely, we have
$$\frac{V(x_i)}{V(x)}=\e^{\int_x^{x_i} S}
\le \e^{C_8 h a^{\alpha-1}}.$$
To keep this ratio bounded, we have to choose $h=O(a^{-\alpha+1}).$ Then we can get the estimate from \eqref{sparkle}.
\begin{equation}
    |{\cal P}f(x)-{\cal P}_{a,h}f(x)| \le
    C V(x_i) [h+ah+\e^{-Ka^\alpha}], \quad \forall x\in [0,a].
\end{equation}
This concludes that $\mcalP_{a,h}\to \mcalP$ as $a\to\infty, h\to 0.$
The authors think that this result is tight since the $O(h)$ error stems mostly from the right-point discretization of the transition kernel, as can be seen from the quantity (I) during the proof of Lemma \ref{start-err}.

With $\mcalP_{a,h}\to\mcalP$ proven, the convergence of invariant measures $\pi_{a,h}\to\pi$ follows naturally from the argument in \S\ref{celestia}. Hence Theorem \ref{main-thm} is completely verified. \qedsymbol

\section{Numerical tests}
In this section, we will perform numerical tests to show the convergence of $\pi_{a,h}$ to $\pi$ of order $O(h)$. The algorithm is as follows.
\begin{enumerate}
    \item Compute the value of $S(x)$ at grid points $x_j$.
    \item Provide any normalized nonnegative vector $\bm\mu\in\BbbR^{N_x/2+1}$ that represents a initial finite volume distribution (sitting in $\ran(I_{a,h})$) on the grid of range $a/2$ and mesh size $h$. In practice we will use a vector with identical positive elements.
    \item Perform the iteration $\bm\mu\leftarrow \bm\mu \widetilde P$ with $\widetilde P$ given by \eqref{ptm} until it converges under the $\ell^1$ norm of vectors.
\end{enumerate}

\begin{example} \label{example1}
Let
\begin{equation}
    g(x)=x, \quad  B(x)=x^2.
\end{equation}
We obtained the results in Figure \ref{numfig}, with varying $a$ and $h$ respectively. It only took around 20 iterations to get to $\pi_{a,h}$, so the convergence of the numerical GF chain to an equilibrium is quite satisfactory.
Note that when plotting $\pi_{a,h}$, we scaled the limiting vector $\bm\pi$ by $h^{-1}$ to convert it to a probability density function.

\begin{figure}[th]
    \centering
    \includegraphics[width=0.45\textwidth]{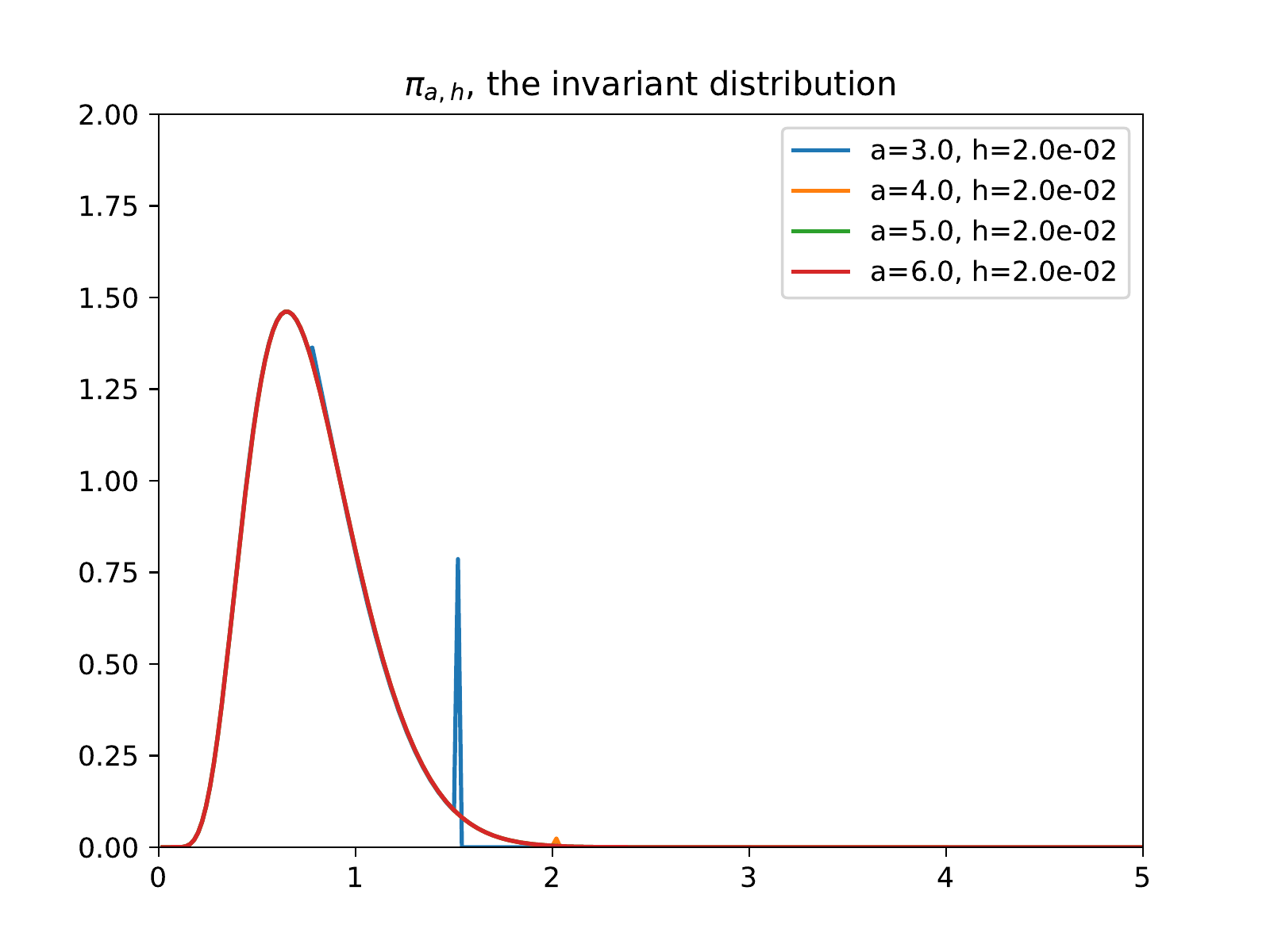}
    \includegraphics[width=0.45\textwidth]{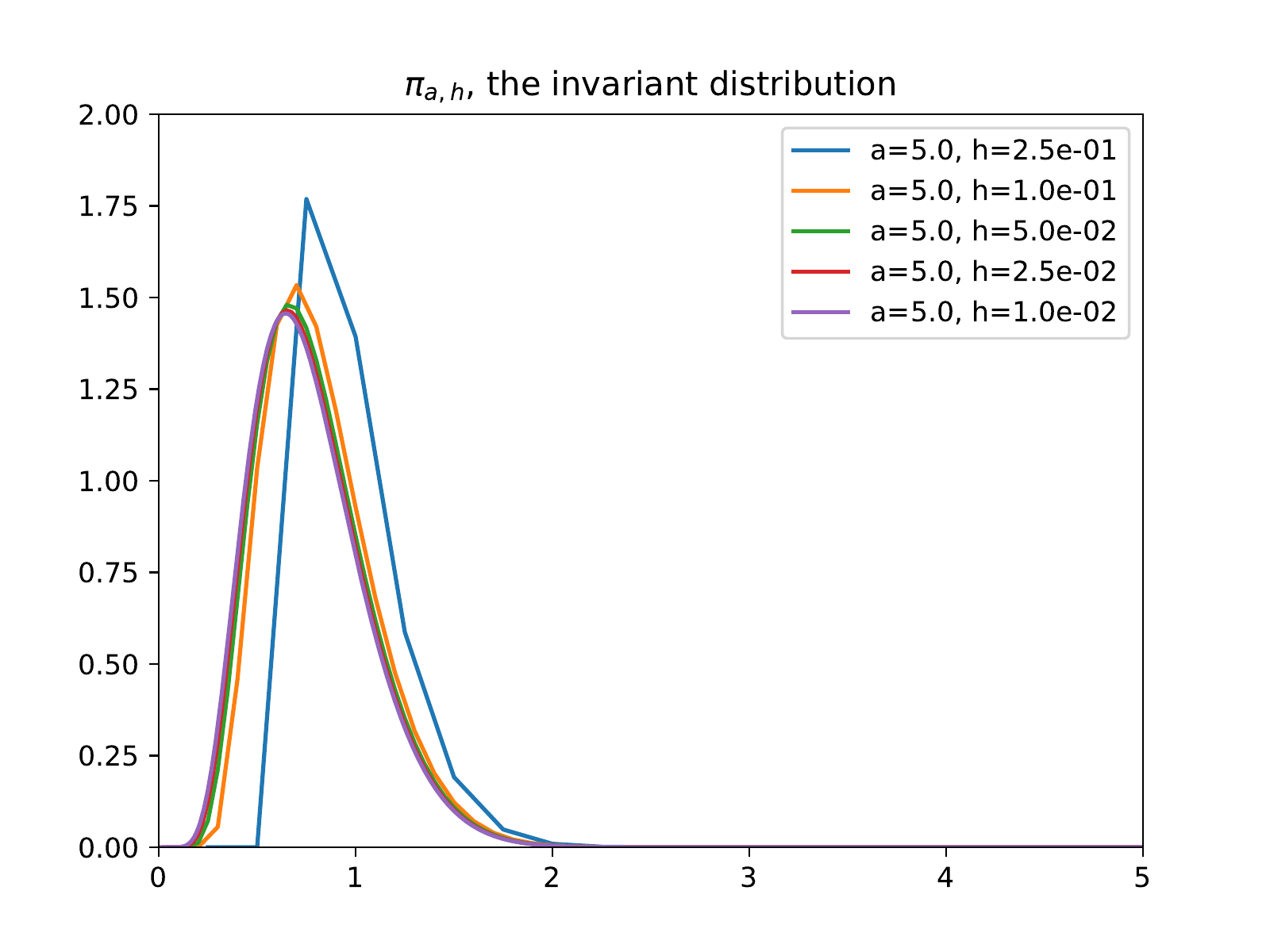}
    \caption{Numerical invariant measure in Example \ref{example1}. Left: different $a$'s and fixed $h=0.02$; right: different $h$'s and fixed $a=5$.} \label{numfig}
\end{figure}

One can observe from Figure \ref{numfig} a quick convergence of $\pi_{a,h}$ to some distribution $\pi$.
In the first figure, the truncation effect seems to have vanished as soon as $a$ reaches 4. Maybe it is because of the rapid decay of the $\e^{-K a^\alpha}$ term in \eqref{main-result}.
In the second figure, the discretization effect also shrinked quickly enough as $h$ got smaller.
Therefore, to obtain a precise approximation of $\nu$, we may first choose $a$ not too large, and then let $h\to 0.$

To quantitatively analyze the errors, we set $a=10$ to make the exponential term negligible and choose different $h$'s. We compute the total-variation error $\|\pi_{a,h}-\pi_{a,2h}\|_{\rm TV}$ to approximate the error $\|\pi_{a,h}-\pi\|$. 
Because the $i$-th unit in the grid of mesh size $h$ occupies half of the $\lfloor\frac{i}{2}\rfloor$-th unit in the grid of mesh size $2h$, we may compute this distance as
$$\|\bm\pi_{a,h}-\bm\pi_{a,2h}\|_1
=\sum_{i=1}^{\frac{N_x}{2}+1} \left| \pi_{a,h}[i] - \frac12 \pi_{a,2h} \left[ \lfloor \nicefrac{i}{2}\rfloor \right] \right|,$$
where $\bm\pi_{a,h}\in\BbbR^{N_x/2+1}, \bm\pi_{a,2h}\in\BbbR^{N_x/4+1}$ are the vector representations of the numerical invariant measures of $\pi_{a,h}$ and $\pi_{a,2h}$ respectively, and out-of-range indexing returns zero.

The result is shown in Table \ref{convtab} and Figure \ref{convfig}. We can clearly see the first order convergence by referring to the reference curve of $y=h$.
\begin{table}[th]
    \centering
    \caption{Numerical error in Example \ref{example1}} \label{convtab}
    \ttfamily
    \begin{tabular}{cccc}
        \toprule
        $h$ &
        1.000e-01 & 5.000e-02 & 2.500e-02 \\
        \midrule
        $\|\pi_{a,h}-\pi_{a,2h}\|_{\rm TV}$ &
        2.0364e-01 & 8.1162e-02 & 3.8011e-02 \\ \midrule
        $h$ & 1.250e-02 & 6.250e-03 & 3.125e-03 \\ \midrule
        $\|\pi_{a,h}-\pi_{a,2h}\|_{\rm TV}$
        & 1.8659e-02 & 9.2620e-03 & 4.6163e-03 \\
        \bottomrule
    \end{tabular}
\end{table} 
\begin{figure}[th]
    \centering
    \includegraphics[width=0.45\textwidth]{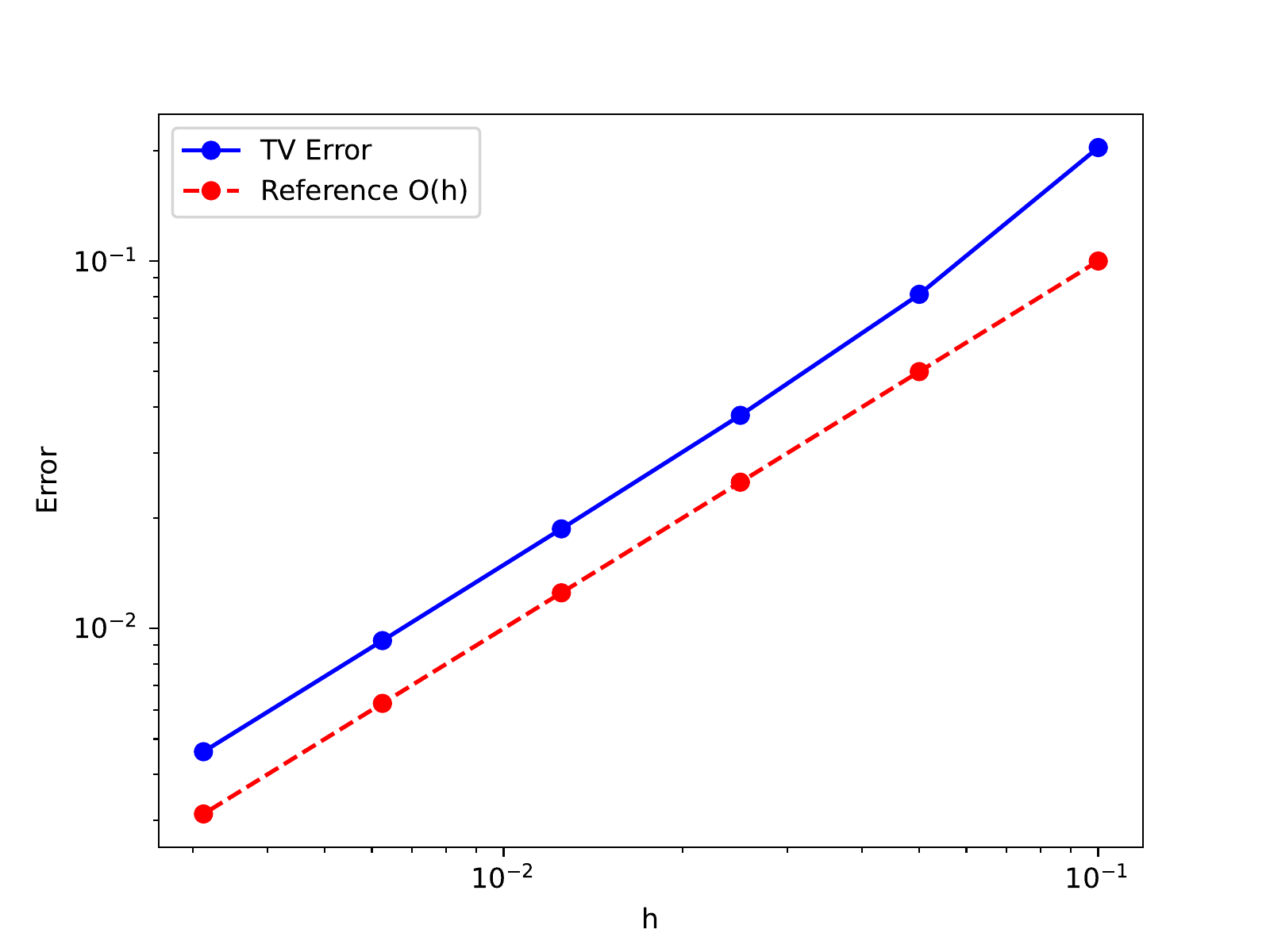}
    \caption{Relation between $\|\pi_{a,h}-\pi_{a,2h}\|_{\rm TV}$ and $h$ in Example \ref{example1} on log-log scale.} \label{convfig}
\end{figure}
\end{example}

We have also tested the algorithm on other choices $B(x)$ and $g(x)$ that do not meet our assumptions \eqref{S-assum} and \eqref{smooth-assum}.

\begin{example} \label{example2}
Letting
$$g(x)=x, \quad B(x)=\max\{x,x^2\},$$
we get a non-smooth but Lipschitz continuous $S(x)$.
Performing the same error analysis as before, we obtain the convergence rate in Figure \ref{figure2}.
\begin{figure}[th]
    \centering
    \includegraphics[width=0.45\textwidth]{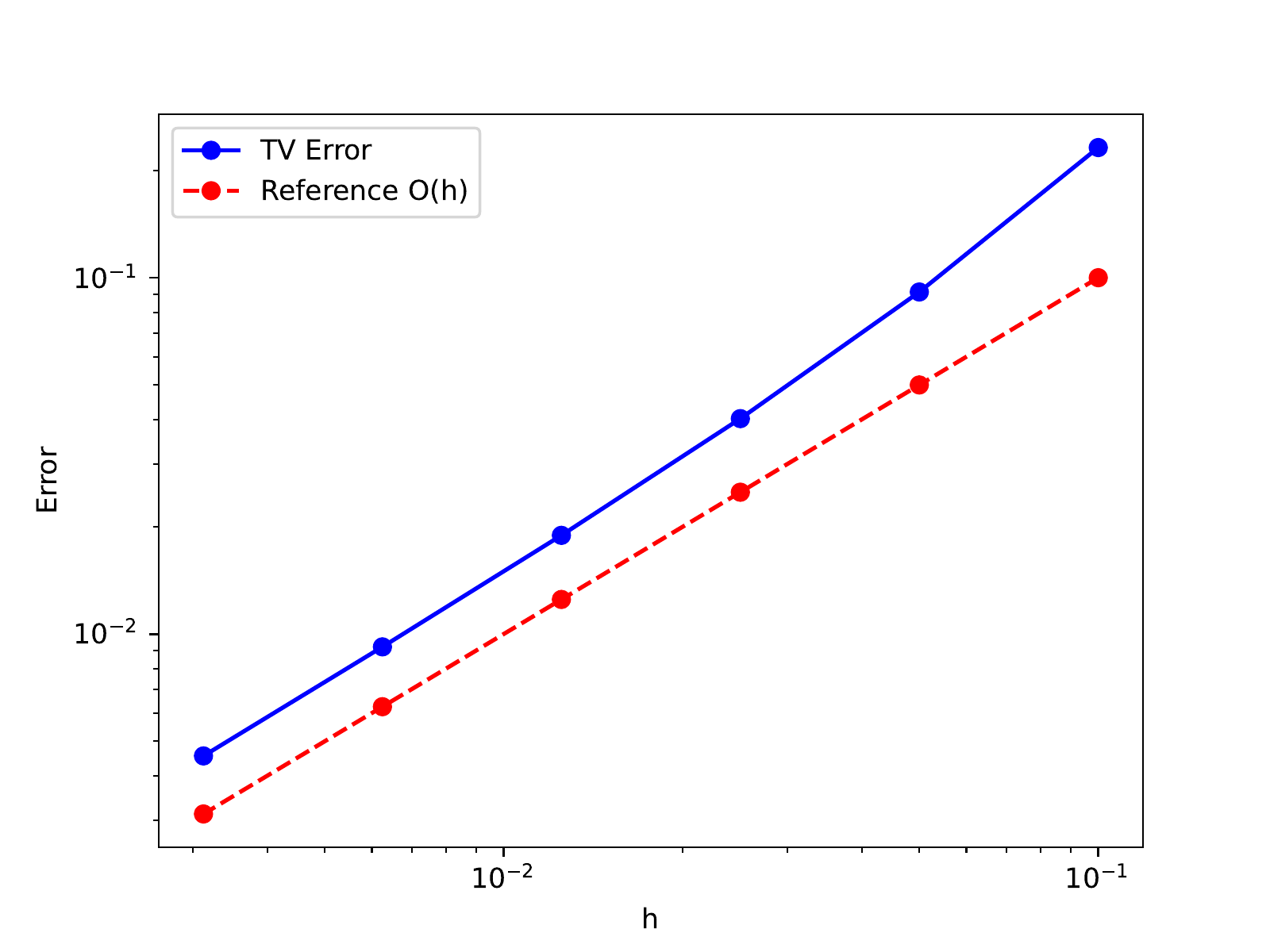}
    \caption{Relation between $\|\pi_{a,h}-\pi_{a,2h}\|_{\rm TV}$ and $h$ in Example \ref{example2}, for Lipschitz continuous $S(x)$.} \label{figure2}
\end{figure}
\end{example}

\begin{example} \label{example3}
Letting
$$g(x)=x, \quad B(x)=x^2+I_{(1,\infty)}(x)$$
we get a discontinuous $S(x)$.
Performing the same error analysis as before, we obtain the convergence rate in Figure \ref{figure3}.
\begin{figure}[H]
    \centering
    \includegraphics[width=0.45\textwidth]{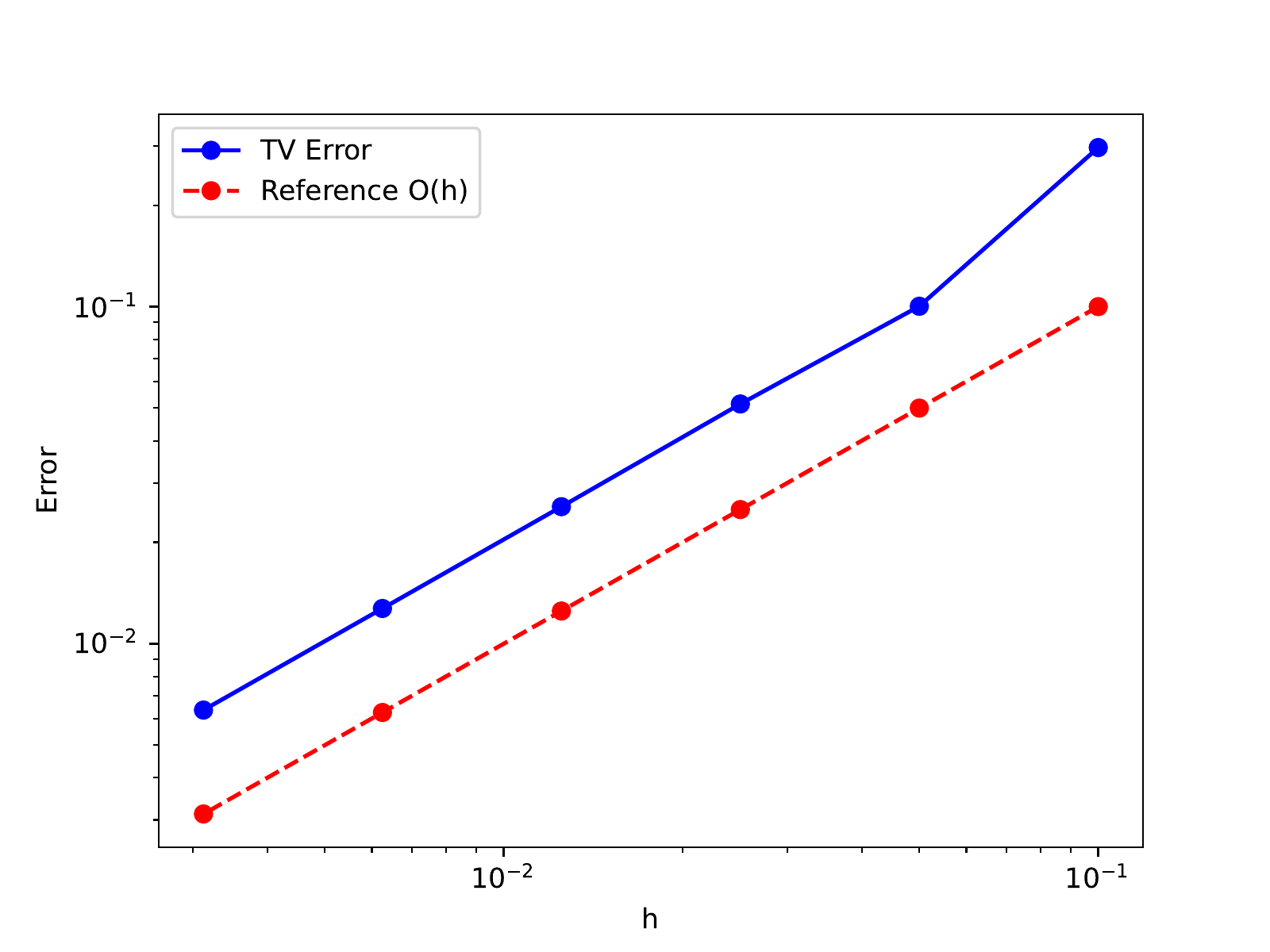}
    \caption{Relation between $\|\pi_{a,h}-\pi_{a,2h}\|_{\rm TV}$ and $h$ in Example \ref{example3}, for discontinuous $S(x)$.} \label{figure3}
\end{figure}
\end{example}

However, the $O(h)$ convergence rate fails when $S(x)$ is singular near $x=0$.

\begin{example} \label{example4}
    Letting
    $$g(x)=x, \quad B(x)=1+x,$$
    which causes $S(x)=1+\frac1x$, a singular function at $x=0$) and running the same program, we get Figure \ref{figure4}.
    There still seems to be a convergence to some $\pi$, but the rate is much slower than $O(h)$. We also find that the convergence of the iteration $\bm\mu\leftarrow\bm\mu\widetilde{P}$ is much slower than before. Our theorem fails in this case.
    \begin{figure}[th]
        \centering
        \includegraphics[width=0.45\textwidth]{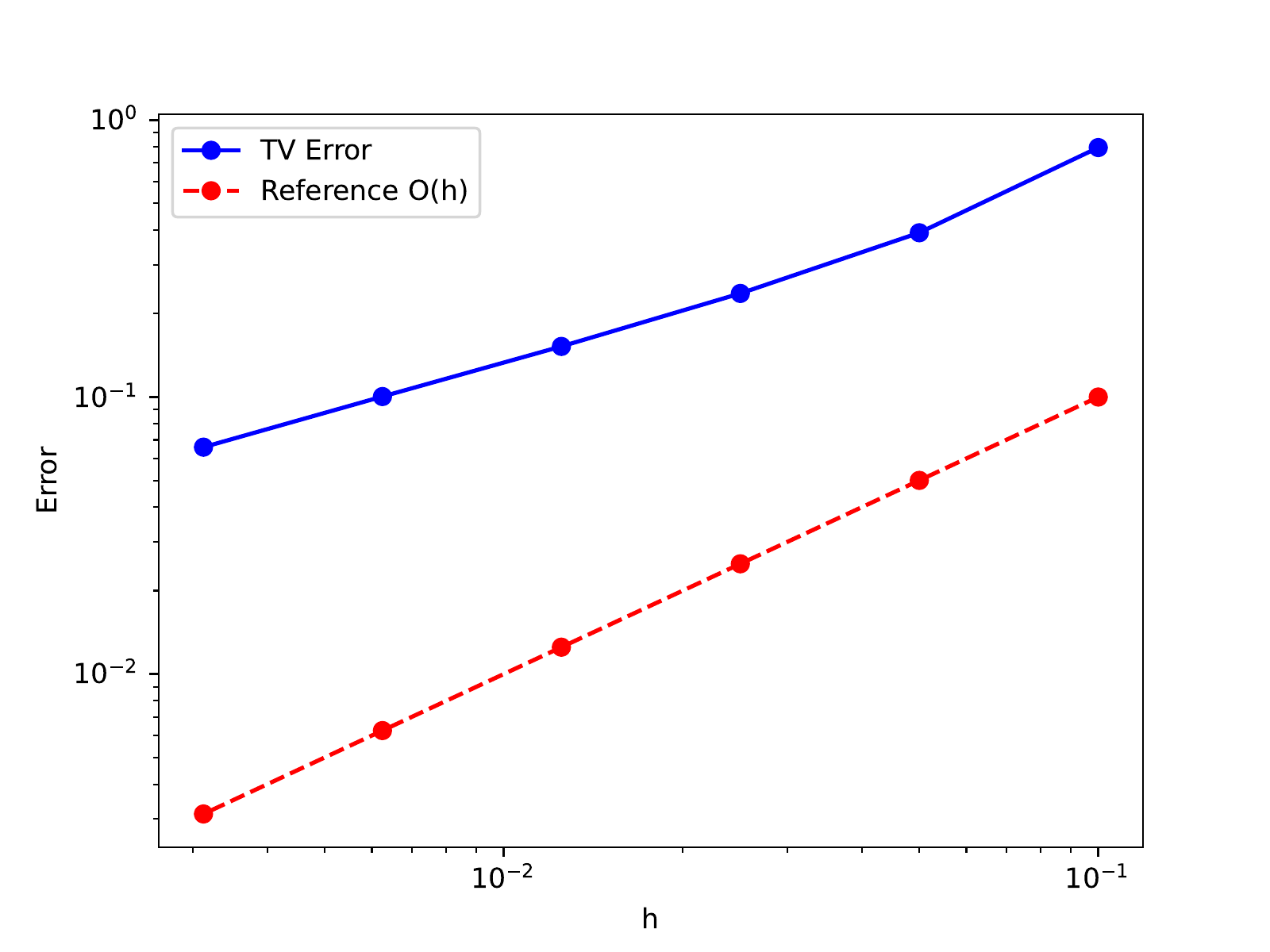}
        \caption{Relation between $\|\pi_{a,h}-\pi_{a,2h}\|_{\rm TV}$ and $h$ in Example \ref{example4}, when choosing singular $S(x)$.}  \label{figure4}
    \end{figure}
\end{example}

\section{Conclusion}
In this paper, we have proposed a justifiable framework under which to generate samples from the forward GF problem. With Theorem \ref{main-thm} proven, we may use it as a generator of test samples to run inverse problem algorithms with, and even actively engage in the solution to inverse problems through Bayesian analysis or maximum likelihood estimate. 
We used a straightforward right-end approximation both in the transition scheme \eqref{num-rule} and in the numerical integration \eqref{Qk}, naturally resulting in an $O(h)$ precision. 

The techniques in our proof are based on the growth rate and smoothness assumptions on the function $S$. The exponentiation of $-S$ in the transitional probability density \eqref{trans-density} provides us with a rapid vanishing rate of the infinite tail, and therefore facilitating much of the order estimate. We note that the smoothness assumption \eqref{smooth-assum} may be weakened by the observations from experiments. 

Future works may involve studies of more general cases such as a singularity of $S(x)$ at $x=0$, as shown in the final numerical example. In that case, our conclusions fail completely such that the $V$-uniform ergodicity of the GF chain and the convergence of the numerical scheme have to be reworked.

\subsection*{Acknowledgements}
This research was partially supported by the National Key R\&D Progrma of China, Project No. 2020YFA0712000 and NSFC Grant No. 12031013, 12171013.

\bibliographystyle{amsplain}
\bibliography{citing.bib}

\end{document}